\documentclass[12pt]{article}

\usepackage{IEEEtrantools}
\usepackage{amsfonts}
\usepackage{amssymb}
\usepackage{amsmath,amsfonts,amsthm}
\usepackage{mathrsfs}
\usepackage{multirow}
\usepackage{url}
\usepackage{tikz}

\newcommand{\bbbt}{\mathbb{T}}

\newcommand{\be}{\begin{equation}}
\newcommand{\ee}{\end{equation}}
\newcommand{\bea}{\begin{eqnarray}}
\newcommand{\eea}{\end{eqnarray}}
\newcommand{\bean}{\begin{eqnarray*}}
\newcommand{\eean}{\end{eqnarray*}}
\newcommand{\brray}{\begin{array}}
\newcommand{\erray}{\end{array}}
\newcommand{\biearray}{\begin{IEEEarray}{rCl}}
\newcommand{\eiearray}{\end{IEEEarray}}

\newcommand{\newsection}[1]{\setcounter{equation}{0}
\setcounter{dfn}{0}
\section{#1}}
\newtheorem{dfn}{Definition}[section]
\newtheorem{thm}[dfn]{Theorem}
\newtheorem{lmma}[dfn]{Lemma}
\newtheorem{ppsn}[dfn]{Proposition}
\newtheorem{crlre}[dfn]{Corollary}
\newtheorem{xmpl}[dfn]{Example}
\newtheorem{rmrk}[dfn]{Remark}

\newcommand{\bdfn}{\begin{dfn}\rm}
\newcommand{\bthm}{\begin{thm}}
\newcommand{\blmma}{\begin{lmma}}
\newcommand{\bppsn}{\begin{ppsn}}
\newcommand{\bcrlre}{\begin{crlre}}
\newcommand{\bxmpl}{\begin{xmpl}}
\newcommand{\brmrk}{\begin{rmrk}\rm}

\newcommand{\edfn}{\end{dfn}}
\newcommand{\ethm}{\end{thm}}
\newcommand{\elmma}{\end{lmma}}
\newcommand{\eppsn}{\end{ppsn}}
\newcommand{\ecrlre}{\end{crlre}}
\newcommand{\exmpl}{\end{xmpl}}
\newcommand{\ermrk}{\end{rmrk}}

\newcommand{\bbc}{\mathbb{C}}

\newcommand{\bbn}{\mathbb{N}}

\newcommand{\cla}{\mathcal{A}}

\newcommand{\clh}{\mathcal{H}}

\newcommand{\cll}{\mathcal{L}}

\newcommand{\prf}{\noindent{\it Proof\/}: }

\def \qed { \mbox{}\hfill
$\Box$\vspace{1ex}}

\begin{document}

%%%%%%%%%%%%%%%%%%%%%%%%%%%%%%%%%

 \author{\sc  { Bipul Saurabh}}
 \title{Quantum Stiefel manifolds}
 \maketitle

%%%%%%%%%%%%%%%%%%%%%%%%%%%%%%%%%%
%%%%%  ABSTRACT
%%%%%%%%%%%%%%%%%%%%%%%%%%%%%%%%%%
 \begin{abstract} 
Quantum analogs of  Stiefel manifolds $SU_{q}(n)/SU_q(n-m)$ were introduced by Podkolzin \& Vainerman. 
The underlying $C^*$-algebra $C(SU_{q}(n)/SU_q(n-m))$   can be described as the $C^*$-subalgebra of $C(SU_q(n))$ 
generated by  elements of last $m$ rows of the fundamental matrix of  $SU_q(n)$. Using $R$-matrix of type $A_{n-1}$, 
one can find certain relations 
involving   elements of last $m$ rows only.  In this paper, by analyzing 
these relations and
using a result of Neshveyev \& Tuset, we establish 
 $C(SU_{q}(n)/SU_q(n-m))$ as a universal $C^*$-algbera given by
 finite sets of generators and relations. 
	
\end{abstract}

{\bf AMS Subject Classification No.:} {\large 58}B{\large 34}, {\large
46}L{\large 87}, {\large
  20}G{\large 42}

\newsection{Introduction}
Compact quantum groups and their homogeneous spaces naturally occur in noncommutative geometry.
A source of 
examples   comes from
$q$-deformation 
of compact simply connected semisimple Lie groups and their quotient spaces.
Many of the $C^*$-algebras arising in this context 
are universal $C^*$-algebras given by finite sets of generators and relations. This fact along with a knowledge of their generators and
relations proved to be very useful in studying topology 
as well as geometry of  these noncommutative spaces. For example, the description 
of the quantum group $SU_q(2)$ and the quantum odd dimensional spheres $S_q^{2n+1}$ as universal $C^*$-algebras
given by  finite sets of generators and relations was used by Vaksman \& Soibelman~(\cite{VakSou-1990ab},\cite{VakSou-1988aa}) 
to compute their $K$-theory. Later using this description, Chakraborty \& Pal (\cite{ChaPal-2003aa},\cite{ChaPal-2008aa}) and 
Pal \& Sundar (\cite{PalSun-2010aa})  constructed spectral triples for these spaces   to study their geometrical aspects. 
Recently 
Saurabh (\cite{Sau-2015aa}) established the quotient space $C(SP_q(2n)/SP_q(2n-2))$ as a universal $C^*$-algebra given by  finite sets of generators 
and relations. Analyzing  the relations among the generators, we get a chain of short exact sequences of $C^*$-algebras  and then using homogeneous extension theory,  we proved in \cite{Sau-2015ab} that 
$C(SP_q(2n)/SP_q(2n-2))$ is isomorphic to $C(S_q^{4n-1})$.

The most well-known examples of quantum homogeneous spaces are quantum Stiefel manifolds $SU_{q}(n)/SU_q(n-m)$. 
Podkolzin \& Vainerman (\cite{PodVai-1999aa}) described all irreducible representations
of the  algebra $S_q^{n,m}$ underlying the manifold $SU_{q}(n)/SU_q(n-m)$. They   proved that the algebra
$S_q^{n,m}$ is generated by matrix entries of last $m$ rows of the fundamental matrix of
$SU_q(n)$. Using $R$-matrix of type $A_{n-1}$, they  found 
certain relations satisfied by these generators of $S_q^{n,m}$ and proved that these relations are full system of relations 
(see Theorem 1, \cite{PodVai-1999aa}). In 
other words, the algebra $S_q^{n,m}$ is the universal algebra generated by $2nm$ number of  generators satisfying 
these relations. 
Using the techniques similar to that used by  Podkolzin \& Vainerman in  \cite{PodVai-1999aa},
one can easily extend some of these results to the $C^*$-algebra level. For example, one can obtain all 
irreducible representations of the $C^*$-algebra underlying the Stiefel manifolds $SU_{q}(n)/SU_q(n-m)$ (in fact Podkolzin \& Vainerman obtained it in  \cite{PodVai-1999aa}).
Similarly  it can be proved 
that $C(SU_{q}(n)/SU_q(n-m))$ is generated by matrix entries of last $m$ rows of the fundamental matrix 
$SU_q(n)$. But extension of  the fact that $S_q^{n,m}$ is the universal algebra generated by elements of last $m$ rows 
satisfying the relations given in Theorem 1, \cite{PodVai-1999aa} to the $C^*$-algebra level is not obvious and demands a proof. 
For $m=1$, 
this was proved by Vaksman \& Soibelman~(\cite{VakSou-1990ab}). 
In this paper, we attempt 
this problem for all $m \leq n$.

Here we should mention that to the best of  our understanding, there seems to be a
 gap in the proof of part (2) of the  Theorem 1 in  \cite{PodVai-1999aa}. Podkolzin \& Vainerman (\cite{PodVai-1999aa}) stated that monomials in 
$\{u_k^l,(u_k^l)^*:1 \leq l \leq n, 1 \leq k \leq n\}$ and $\{u_k^l,(u_k^l)^*:n-m+1 \leq l \leq n, 1 \leq k \leq n\}$ having 
the lexicographic order form  bases of the algebras  $U_q(n)$ underlying the manifold $SU_q(n)$ and 
$S_q^{n,m}$ respectively.  But since $\sum_{k=1}^nu_k^l(u_k^j)^*=\delta_{lj}$, this is not true. In fact it follows from 
diamond lemma (see page 103, \cite{KliSch-1997aa}) that  reduced monomials in  the generators of $U_q(n)$ and 
$S_q^{n,m}$  will 
form bases of the algebras $U_q(n)$ and 
$S_q^{n,m}$ respectively. But then it is not obvious that a reduced monomial  
of the generators $\{u_k^l,(u_k^l)^*:n-m+1 \leq l \leq n, 1 \leq k \leq n\}$ in the algebra  $S_q^{n,m}$ will be 
a reduced monomial  in the algebra  $U_q(n)$ as there are more relations among the generators of the algebra  $U_q(n)$. One thing we should point out here that Podkolzin \& Vainerman (\cite{PodVai-1999aa}) used  left coset space 
and as a consequence,  elements of 
last $m$ columns of the fundamental matrix of $SU_q(n)$ are the generators of the underlying algebra $S_q^{n,m}$. Here we take right coset space  $SU_{q}(n)/SU_q(n-m)$ and hence elements of 
last $m$ rows are the generators of the underlying algebra $S_q^{n,m}$. But similar arguments hold in both cases.

To tackle the problem, we first write down all irreducible representations of the $C^*$-algebra $C(SU_{q}(n)/SU_q(n-m))$ by 
applying a result of  Neshveyev \& Tuset (\cite{NesTus-2012ab}). We associate some diagrams to each irreducible representation and 
using this, we describe certain properties of these representations. We then define the  $C^*$-algebra $C^{n,m}$ as a universal $C^*$-algebra
generated by $nm$ generators satisfying those relations satisfied by the generators of $C(SU_{q}(n)/SU_q(n-m))$.  Due to
universal property of $C^{n,m}$, we get a surjective homomorphism from $C^{n,m}$ to $C(SU_{q}(n)/SU_q(n-m))$.
 By analyzing the relations carefully and using some properties of irreducible representations of $C(SU_{q}(n)/SU_q(n-m))$,  we show that  all
irreducible representations of the $C^*$-algebra $C^{n,m}$  factor through this homomorphism. This implies that the homomorphism
is injective. Therefore, the two $C^*$-algebras $C^{n,m}$  and $C(SU_{q}(n)/SU_q(n-m))$ are isomorphic which establishes $C(SU_{q}(n)/SU_q(n-m))$ 
as a universal $C^*$-algebra given by  finite sets of generators and relations.

We set up some notations which will be used throughout this paper. The standard basis of the 
Hilbert space $L_2(\bbn)$ will be denoted by
$\left\{e_n: n\in \bbn \right\}$. 
We  denote the left shift operator and the number operator on $L_2(\bbn)$ by $S$ and $N$ respectively.  Let $p$ denote the rank one projection sending 
$e_0$ to $e_0$. We denote by 
$\bbbt^n$ the $n$-dimensional torus. For $\alpha= (\alpha_1, \cdots ,\alpha_n) \in \bbn^n$, the symbol $|\alpha |$ denotes the number
$\sum_{i=1}^n \alpha_i$. The number 
 $q$ will denote a real number in the interval  $(0,1)$.

 \newsection{ODQS(q) relations}
 In this section, we describe  some properties of  ordered tuples of bounded operators satisfying  certain relations. 
 Let $T_1,T_2,\cdots T_n$ be bounded linear  operators acting on  
 a Hilbert space $\clh$. We say that the ordered tuple  $(T_1,T_2, \cdots ,T_n)$ satisfies odd dimension quantum sphere relations with parameter $q$
  abbreviated as $ODQS(q)$ relations  if it obeys the following relations:
 \begin{align*}
	T_iT_j &= qT_jT_i,                 &    &  1\leq j < i\leq n,  \\
	T_i^{*}T_j &= qT_jT_i^{*},  &  & 1\leq i \neq j\leq n,  \\
	T_i^{*}T_i - T_iT_i^{*} &=  (1-q^2)\sum_{j>i}T_jT_j^{*},   
	      & \sum_{i=1}^{n}T_iT_i^{*}&=1. 
\end{align*}
It follows from  above relations that there exists a unique number $\ell$ in $\{1,2, \cdots, n\}$ such that $T_{\ell}$ is a nonzero normal operator and for 
$i>\ell$, $T_i=0$. We call this number $\ell$ the rank of $(T_1,T_2, \cdots ,T_n)$. 
We say that the tuple $(T_1, \cdots, T_n)$  is irreducible if there is no proper subspace of $\clh$  invariant 
under the action of $\{T_i, T_i^*: i \in \{1,2,\cdots,n\}\}$. 

\bppsn \label{ppsn-spectrum}
Suppose that the tuple  $(T_1,T_2, \cdots ,T_n)$ obeys $ODQS(q)$ relations and of rank $\ell$. Let $\omega$ be the operator $T_{\ell}^*T_{\ell}$ 
and  $\clh_0$ be the eigenspace of $\omega$ corresponding to the eigenvalue $1$.  Then one has 
\begin{enumerate}
\item for $1\leq i <\ell$, \quad 
$ T_i\omega=q^{-2}\omega T_i \qquad $ and $ \qquad T_i^*\omega =q^2\omega T_i^*$,

\item 
$\omega = I \qquad \mbox{ on } \;\bigcap_{i=1}^{\ell-1} \ker T_{i}^{*}$,
\item 
$1_{(q^{2m+2}, q^{2m})}(\omega) = 0 \quad \forall m \in \bbn$,
\item
$\ker T_{i}\subseteq \ker T_{k}^{*} \mbox{ for } k\geq i \mbox{ and } 1 \leq i \leq \ell$,
\item
if $u$ is a nonzero eigenvector of $\omega$ corresponding to the eigenvalue $q^{2m}$, then 
$u \notin \ker T_{i}$  for  $1 \leq i \leq \ell-1$.
\item
$\sigma (\omega) = \left\{q^{2m}: m \in \bbn \right\} \bigcup \left\{0\right\}$. 
\item
For $i \in \{1,2,\cdots,\ell-1\}$ and $m \in \bbn-\{0\}$, one has
\[
 T_i^*T_i^{m}=T_i^{m}T_i^*+(1-q^{2m})\sum_{j>i}T_i^{m-1}T_jT_j^*.
\]

\item 
 Let
$ \clh_{(\alpha_1,\cdots \alpha_{\ell-1})}:=T_1^{\alpha_1}T_2^{\alpha_2} \cdots T_{\ell-1}^{\alpha_{\ell-1}}\clh_0$
where $\alpha_i \in \bbn $ for all $i \in \{1,2,\cdots \ell-1\}$.
Then for different values of $\alpha=(\alpha_1,\cdots \alpha_{\ell-1})$, the subspaces $\clh_{(\alpha_1,\cdots ,\alpha_{\ell-1})}$ are  nonzero and mutually 
orthogonal. 
\item
If in addition,  $\ker(\omega)=\{0\}$, then one has
\[
\clh=\oplus_{\alpha_i \in \bbn}\clh_{(\alpha_1,\cdots \alpha_{\ell-1})}=
 \oplus_{\alpha_i \in \bbn}T_1^{\alpha_1}T_2^{\alpha_2} \cdots T_{\ell-1}^{\alpha_{\ell-1}}\clh_0.
\]
 \item
 The tuple $(T_1, \cdots, T_n)$ is irreducible if and only if $\clh_0$ is one dimensional. Moreover in this case,  $\ker(\omega)=0$ and for a nonzero vector $h \in \clh_0$, 
 \[\Big\{\frac{T_1^{\alpha_1}T_2^{\alpha_2} \cdots T_{\ell-1}^{\alpha_{\ell-1}}h}{\|T_1^{\alpha_1}T_2^{\alpha_2} \cdots T_{\ell-1}^{\alpha_{\ell-1}}h\|}: \alpha_i \in \bbn \mbox{ for all } i \in \{1,2,\cdots ,\ell-1\} \Big\}\]
form an orthonormal basis of $\clh$.
Further for $h \in \clh_0$,  $T_{\ell}h=th$ for some $t \in \bbbt$. 
We call $t \in \bbbt$ the angle of the irreducible tuple $(T_1, \cdots, T_n)$.
\item
For $1\leq i <\ell$ and $h \in \clh_0$, one has 
\begin{displaymath}
 T_i^*h=0, \qquad \mbox{ and } \qquad  T_i^*T_i=(1-q^2)h.
\end{displaymath}

\end{enumerate}
\eppsn 
\prf
\begin{enumerate}
\item
Easy to see from ODQS(q) relations.
 \item 
 It follows from the relation $\sum_{i=1}^{\ell}T_iT_i^{*}=1$.
 \item
 From the ODQS(q) relations, it follows that  
$T_i^*f(\omega)=f(q^2\omega)T_i^*$ for all $i \in \{1,2, \cdots,\ell-1\}$ for all
continuous functions $f$ and hence for all $L_\infty$ functions.
Thus 
\begin{IEEEeqnarray}{rCl}
         T_i^*1_{(q^{2n+2},q^{2n})}(\omega) &=&1_{(q^{2n+2},q^{2n})}(q^2\omega)T_i^* \nonumber \\
                       &=&1_{(q^{2n},q^{2n-2})}(\omega)T_i^*. \nonumber
\end{IEEEeqnarray}
By repeated application and using the relation $\sum_{i=1}^{\ell}T_iT_i^{*}=1$ and the fact that
$\sigma(\omega)\subseteq[0,1]$, it follows that $1_{(q^{2n+2},q^{2n})}(\omega)=0$.
\item
Let $h \in \ker(T_i)$. Then we have
\[
\left\langle T_i^*T_ih, h \right\rangle = \left\langle T_iT_i^*h + (1-q^2)\sum_{k>i}T_kT_k^*h, h \right\rangle.
\]
Hence
\[
\left\|T_i^*h\right\|^2 + (1-q^2)\sum_{k>i}\left\|T_k^*h\right\|^2= 0.
\]
Hence $\left\|T_k^*h\right\| = 0$ for all  $k\geq i$, which means 
 $h  \in  \ker(T_{k}^{*})$ for all  $k\geq i$. 

\item
From part $(4)$, we have $\ker(T_{i})\subseteq \ker(T_{\ell}^*) = \ker(T_{\ell})= \ker(\omega)$.
Now if $u$ is a non-zero eigenvector of $\omega$ corresponding to eigenvalue $q^{2m}$ for
some $m \in \bbn$, then $u \notin \ker(T_{\ell}^*)$. Hence $u \notin \ker(T_{i})$ for $1 \leq i \leq \ell$. 
\item
From part $(3)$ and the fact that $\left\|\omega\right\|\leq 1$, it follows that the spectrum $\sigma(\omega)$ of $\omega$ is contained in 
$\left\{q^{2m}: m \in \bbn\right\} \bigcup \left\{0\right\}$. Define
\[
 A = \left\{m \in \bbn : q^{2m} \in  \sigma(\pi(\omega))\right\}.
\]
Since $\omega \neq 0$,  $A \neq \emptyset$. Let $m_{0}=\inf\{m \in \bbn : q^{2m} \in  \sigma(\omega)\}$.
Let $u$ be a nonzero eigenvector corresponding to $q^{2m_{0}}$. Assume $u \notin  \ker(T_{i}^{*})$
for some $i \in \left\{1,2,\cdots ,\ell-1\right\}$. Then from part (1), it follows that  $T_{i}^{*}u$ 
is a nonzero eigenvector corresponding to the eigenvalue $q^{2m_{0}-2}$, which contradicts the 
fact that $m_{0}$ is $\inf A$. Hence $u \in \bigcap_{i=1}^{\ell-1}\ker(T_{i}^{*})$.
As  $\omega = I$ on $\bigcap_{i=1}^{\ell-1}\ker(T_{i}^{*})$, we get $m_{0} = 0$. From part $(5)$, it follows that 
$u \notin \ker(T_{i})$ for any $i \in \left\{1,2,\cdots \ell\right\}$. Hence we have
$T_1^{m}u$ is a nonzero eigenvector corresponding to eigenvalue $q^{2m} \mbox{ for all } m \in \bbn$. This proves the claim. 
\item
It follows by applying  the relation $T_i^{*}T_i = T_iT_i^{*} +  (1-q^2)\sum_{j>i}T_jT_j^{*}$ repeatedly.
\item 
For any nonzero $h \in \clh_0$, the vector $T_1^{\alpha_1}\cdots T_{\ell-1}^{\alpha_{\ell-1}}h$ is an eigenvector corresponding to eigenvalue 
$q^{2|\alpha|}$. Hence it follows from part $(5)$ that $\clh_{(\alpha_1,\cdots \alpha_{\ell-1})}$ is nonzero. Let 
$\alpha=(\alpha_1,\cdots,\alpha_{\ell-1})$ and $\alpha^{'}=(\alpha_1^{'},\cdots,\alpha_{\ell-1}^{'})$ be two different tuple of 
positive integers. Without loss of generality, we assume that $\alpha_1 \neq 0$. For $h,h^{'} \in \clh_0$, we have
\begin{IEEEeqnarray}{lCl}
 \langle T_1^{\alpha_1} \cdots T_{\ell-1}^{\alpha_{\ell-1}}h,  T_1^{\alpha_1^{'}}\cdots T_{\ell-1}^{\alpha_{\ell-1}^{'}}h^{'} \rangle \nonumber \\
 =\langle T_1^{\alpha_1-1} \cdots T_{\ell-1}^{\alpha_{\ell-1}}h,  T_1^*T_1^{\alpha_1^{'}}\cdots T_{\ell-1}^{\alpha_{\ell-1}^{'}}h^{'} \rangle \nonumber \\
 =q^{\sum_{i=2}^{\ell-1}\alpha_i^{'}}(1-q^{2\alpha^{'}})
 \langle T_1^{\alpha_1-1} \cdots T_{\ell-1}^{\alpha_{\ell-1}}h, T_1^{\alpha_1^{'}-1}\cdots T_{\ell-1}^{\alpha_{\ell-1}^{'}}h^{'} \rangle 
  \nonumber
\end{IEEEeqnarray}
Last equality follows from part (7) and the ODQS(q) relations. Now by induction on $|\alpha|$, we get the claim.

\item
Let $\clh_m$ be the eigenspace of $\omega$ corresponding to the eigenvalue $q^{2m}$. Then by spectral theorem and the fact that $\ker \omega=\{0\}$, one has
 $\clh=\oplus_{m \in \bbn}\clh_m.$
So, to prove the claim, we need to show that for all $m \in \bbn$, 
\[
\clh_m=\oplus_{\{\alpha_i \in \bbn:\sum_{i=1}^{\ell-1}\alpha_i=m\}}T_1^{\alpha_1}T_2^{\alpha_2} \cdots T_{\ell-1}^{\alpha_{\ell-1}}\clh_0.
\] 
From the ODQS(q) relations, we have  $\oplus_{\{\alpha_i \in \bbn:\sum_{i=1}^{\ell-1}\alpha_i=m\}}T_1^{\alpha_1}T_2^{\alpha_2} \cdots T_{\ell-1}^{\alpha_{\ell-1}}\clh_0 \subset \clh_m$. Take 
 a   nonzero vector $h$ in $\clh_m$ for $m \neq 0$. 
 Then  from the relation $\sum_{i-1}^{\ell}T_iT_i^*=1$, it follows that there exists $i \in \{1,2, \cdots, \ell-1\}$ such that $T_i^*h \neq 0$. Hence we get a nonzero vector $T_i^*h$ in 
$\clh_{m-1}$. Proceeding in this way and arranging term using ODQS(q) relations, we get $\beta_i \in \bbn$ such that $\sum_{i=1}^{\ell-1}\beta_i=m$ and $h_0:=(T_{\ell-1}^*)^{\beta_{\ell-1}} \cdots (T_1^*)^{\beta_1}h$ 
is nonzero vector in $\clh_0$. Therefore, we get
\begin{IEEEeqnarray}{rCl}
 \langle h, T_1^{\beta_1}T_2^{\beta_2} \cdots T_{\ell-1}^{\beta_{\ell-1}}h_0 \rangle &=&  \langle  (T_{\ell-1}^*)^{\beta_{\ell-1}} \cdots (T_1^*)^{\beta_1}h, h_0)\rangle \nonumber \\
 &=& \langle h_0, h_0 \rangle \nonumber \\
 &\neq & 0. \nonumber
\end{IEEEeqnarray} 
This implies that  no  nonzero vector in $\clh_m$ is orthogonal to the vector subspace   $\oplus_{\{\alpha_i \in \bbn:\sum_{i=1}^{\ell-1}\alpha_i=m\}}T_1^{\alpha_1}T_2^{\alpha_2} \cdots T_{\ell-1}^{\alpha_{\ell-1}}\clh_0$. 
Hence we get 
\[
 \clh_m \subset \oplus_{\{\alpha_i \in \bbn:\sum_{i=1}^{\ell-1}\alpha_i=m\}}T_1^{\alpha_1}T_2^{\alpha_2} \cdots T_{l-1}^{\alpha_{\ell-1}}\clh_0.
\]
This proves the claim.

\item
 It follows from the relations $\sum_{i=1}^{\ell}T_iT_i^*=1$ that $T_l$ acts as a unitary operator on $\clh_0$. Hence there exists an orthonormal basis $\{h_i\}_{i \in \bbn}$ of $\clh_0$
such that $T_{\ell}h_i=t_ih_i$ for 
some $t_i \in \bbbt$. It is easy to see that for any $i \in \bbn$,  the subspace spanned by the nonzero vectors 
$\{T_1^{\alpha_1}T_2^{\alpha_2} \cdots T_{\ell-1}^{\alpha_{\ell-1}}h_i:\alpha_i \in \bbn\}$ is an invariant subspace. This proves that if $(T_1, \cdots, T_n)$ is irreducible then 
$\clh_0$ is one dimensional. To show the other way, take $h \in \clh$. As done in part (9), one can show that by applying $(T_i)^*$'s repeatedly  for appropriately chosen 
$ i \in \{1,2,\cdots ,\ell-1\}$, one can take $h$ to $\clh_0$. 
This implies that any nonzero invariant subspace contains a vector in $\clh_0$. Since $\clh_0$ is one dimensional, any nonzero invariant subspace contains $\clh_0$. This implies that  the only nonzero invariant subspace is $\clh$ and
the tuple  $(T_1,T_2, \cdots, T_n)$ is irreducible.\\
From the ODQS($Q$) relations, one can easily see that $\ker(\omega)$ is an invariant subspace. Since $\omega \neq 0$, it follows that 
$\ker(\omega)=\{0\}$. Other parts  of the claim follow from part (8) and part (9).  

\item 
 It follows from ODQS(q) relations that for $1 \leq i < \ell$ and $h \in \clh_0$, $\omega (T_i^{*}h)=\frac{1}{q}T_i^{*}h$.
Hence from part (6), we get $T_i^*h=0$. Using this fact and the relation $T_i^{*}T_i = T_iT_i^{*} +  (1-q^2)\sum_{j>i}T_jT_j^{*}$, we get  $T_i^*T_ih=(1-q^2)h$. 

\end{enumerate}
\qed 
\brmrk
Note that the operators $T_i$'s occuring in the tuple $(T_1, \cdots, T_n)$ that satisfies ODQS(q) relations are  homomorphic image of the generators 
of the $C^*$-algebra $C(S_q^{2n-1})$ of continuous functions on  odd dimensional quantum sphere.  Using this fact along with the knowledge of irreducible representations of   the $C^*$-algebra $C(S_q^{2n-1})$, 
one can prove above proposition but to make the paper self-contained, we are giving a direct proof.  
\ermrk

We shall now define some notions related to the  ordered $n$-tuples satisfying ODQS(q) relations. Let $\{(T_1^i,T_2^i,\cdots,T_n^i)\}_{i \in I}$ 
be ordered $n$-tuples obeying ODQS(q) relations with rank $\ell_i$. Their direct sum $\oplus_{i \in I}(T_1^i,T_2^i,\cdots,T_n^i)$ is 
defined as 
$(\oplus_{i \in I}T_1^i,\cdots,\oplus_{i \in I}T_n^i)$. Clearly rank of the direct sum of ordered $n$-tuples  is $\max\{\ell_i:i \in I\}$. 
Moreover, one can show that  if $\ker(\oplus_{i \in I}T_{l}^i)=\{0\}$ where $\ell:=\max\{l_i:i \in I\}$ then for all $i \in I$, $\ell_i=\ell$. 
We say that two ordered $n$-tuple $(T_1,T_2,\cdots,T_n)$ and $(T_1^{'},T_2^{'},\cdots,T_n^{'})$ of operators acting on 
the Hilbert spaces $\clh$ and $\clh^{'}$ respectively that satisfy  
ODQS(q) relations are isomorphic if there exists unitary $U:\clh \rightarrow \clh^{'}$ such that $UT_iU^*=T_i^{'}$ for 
all $1\leq i \leq n$. The following proposition says that isomorphism class of an irreducible ordered $n$-tuple is characterized 
by its rank and angle. More precisely, 
\bppsn \label{isomordered}
Let $(T_1,T_2,\cdots,T_n)$ and $(T_1^{'},T_2^{'},\cdots,T_n^{'})$ be two irreducible ordered $n$-tuple obeying ODQS(q) relations  on 
the Hilbert spaces $\clh$ and $\clh^{'}$ respectively. Let their    ranks be $\ell$ and $\ell^{'}$ respectively and angles be  $t$ and $t^{'}$ respectively. 
Then $(T_1,T_2,\cdots,T_n)$ and $(T_1^{'},T_2^{'},\cdots,T_n^{'})$ are isomorphic if and only if $\ell=\ell^{'}$ and $t=t^{'}$.
\eppsn
\prf It is easy to show that if $(T_1,T_2,\cdots,T_n)$ and $(T_1^{'},T_2^{'},\cdots,T_n^{'})$ are isomorphic then  $\ell=\ell^{'}$ and $t=t^{'}$. To show the other way, let
$\clh_0$  and $\clh_0^{'}$ be the eigenspaces  of $T_{\ell}^*T_{\ell}$ and $(T_{\ell}^{'})^*T_{\ell}^{'}$ respectively corresponding to the eigenvalue $1$. From 
 part (10) of the proposition \ref{ppsn-spectrum}, $\clh_0$  and $\clh_0^{'}$ are one dimensional. Let $h$ and $h^{'}$ be  unit vectors in 
 $\clh_0$  and $\clh_0^{'}$ respectively. Let
 \[
  h_{(\alpha_1,\cdots \alpha_{\ell-1})}:=
  \frac{T_1^{\alpha_1}T_2^{\alpha_2} \cdots T_{\ell-1}^{\alpha_{\ell-1}}h}{\|T_1^{\alpha_1}T_2^{\alpha_2} \cdots T_{\ell-1}^{\alpha_{\ell-1}}h\|}, \quad \qquad 
  h_{(\alpha_1,\cdots \alpha_{\ell-1})}^{'}:=
 \frac{(T_1^{'})^{\alpha_1}(T_2^{'})^{\alpha_2} \cdots (T_{\ell-1}^{'})^{\alpha_{\ell-1}}h^{'}}{\|(T_1^{'})^{\alpha_1}(T_2^{'})^{\alpha_2} \cdots (T_{\ell-1}^{'})^{\alpha_{\ell-1}}h^{'}\|}
 \]
It follows from part (10) of proposition \ref{ppsn-spectrum}  that 
$\{h_{(\alpha_1,\cdots \alpha_{\ell-1})}: \alpha_i \in \bbn \mbox{ for } 1\leq i\leq \ell-1\}$ and 
$\{h_{(\alpha_1,\cdots \alpha_{\ell-1})}^{'}: \alpha_i \in \bbn \mbox{ for } 1\leq i\leq \ell-1\}$ are 
orthonormal bases of $\clh$ and $\clh^{'}$ respectively. Let  
$U:\clh \rightarrow \clh^{'}$ such that $Uh_{(\alpha_1,\cdots \alpha_{\ell-1})}=h_{(\alpha_1,\cdots \alpha_{\ell-1})}^{'}$. Using 
ODQS$(q)$ relations and the fact that $t=t^{'}$, it is easy to verify that $UT_iU^*=T_i^{'}$ for 
all $1\leq i \leq n$. This completes the proof. 
\qed 

\brmrk \label{torusinvariant}
Suppose that  $(T_1, \cdots, T_n)$ obeys ODQS(q) relations on $\clh$  and of rank $\ell$. Then  it is easy to see that 
the tuple $(\omega_1T_1,\omega_2T_2,\cdots , \omega_nT_n)$ also satisfies ODQS(q) relations with rank $\ell$ where $\omega_i\in \bbbt$ for all $1\leq i \leq n$. Further, it follows 
from  proposition \ref{isomordered} that the two tuples $(\omega_1T_1,\omega_2T_2,\cdots , \omega_nT_n)$ and 
$(\omega_1^{'}T_1,\omega_2^{'}T_2,\cdots , \omega_n^{'}T_n)$ are isomorphic if and only if $\omega_{\ell}=\omega_{\ell}^{'}$.
\ermrk

\newsection{Representation theory of $C(SU_{q}(n)/SU_q(n-m))$}
In the present section, we recall  representation theory  of $C(SU_{q}(n)/SU_q(n-m))$ using a result of Neshveyev \& Tuset and then
associate a diagram to each irreducible representation. 
Let the fundamental matrices of the compact quantum groups  $SU_q(n)$ and $SU_q(n-m)$ be  $U:=(\!(u_{k}^{l})\!)$ and 
$V:=(\!(v_{k}^{l})\!)$ respectively. Consider the following map:
\[
 \phi(u_k^l)=\begin{cases}
               v_k^l & \mbox{ if } 1\leq l,k \leq n-m, \cr
               1 &  \mbox{ if } n-m+1 \leq l=k \leq n, \cr
               0 & \mbox{ otherwise. }
              \end{cases}
\]
One can check that $\phi$ is a quantum group homomorphism from $SU_q(n)$ onto $SU_q(n-m)$ and hence  $SU_q(n-m)$ is a subgroup of $SU_q(n)$. 
The $C^*$-algebra 
of continuous functions on the right coset space $SU_q(n)/SU_q(n-m)$  is given by
\begin{IEEEeqnarray}{rCl}
 C(SU_q(n)/SU_q(n-m))=\left\{a \in C(SU(n)) : (\phi \otimes id)\Delta(a)=I \otimes a\right\} \label{right quotient}
\end{IEEEeqnarray}
where $\Delta$ is the comultiplication map of $SU(n)$. The following theorem gives finite number of  
generators of the quotient space $C(SU_q(n)/SU_q(n-m))$. It  can be considered as an extension of 
Theorem 1, \cite{PodVai-1999aa} from algebra level to $C^*$-algebra level. One can  easily derive it from the  proof 
of Theorem 1, \cite{PodVai-1999aa} but for the sake of
completeness, we are giving its proof here.
\bthm \cite{PodVai-1999aa}\label{quotient}
The quotient space $C(SU_{q}(n)/SU_{q}(n-m))$ is the $C^*$-algebra generated by 
$\left\{u_k^l : n-m+1 \leq l \leq n, 1 \leq k \leq n \right\}$.
\ethm 
\prf 
Using the definition of quotient space, one can easily check that  
%$\left\{u_k^l : l \in \{1,2,\cdots m\}, k \in \{1,2,\cdots n\}\right\}$
$u_k^l$'s are in $C(SU_{q}(n)/SU_{q}(n-m))$ for $n-m+1 \leq l \leq n$ and  $1 \leq k \leq n $. So, 
\[
C(SU_{q}(n)/SU_{q}(n-m)) \supseteq C^*\left\{u_k^l : n-m+1 \leq l \leq n, 1 \leq k \leq n \right\}.
\]
 To show the equality,
consider the co-multiplication action on $C(SU_{q}(n)/SU_{q}(n-m))$ by the compact quantum group $C(SU_q(n))$ given by
\begin{IEEEeqnarray}{rCl}
 C(SU_{q}(n)/SU_{q}(n-m)) &\longrightarrow & C(SU_{q}(n)/SU_{q}(n-m)) \otimes C(SU_q(n)) \nonumber \\
 a &\longmapsto & \Delta a. \nonumber 
\end{IEEEeqnarray}
By theorem $1.5$,  \cite{Pod-1995aa}, we get
\begin{IEEEeqnarray}{rCl}
C(SU_{q}(n)/SU_{q}(n-m))= \overline{\oplus_{\lambda \in \widehat{SU(n)}}\oplus_{i \in I_{\lambda} }W_{\lambda,i}} \nonumber 
\end{IEEEeqnarray}
where $\lambda$ represents a finite-dimensional irreducible co-representation $u^{\lambda}$ of $C(SU_q(n))$, $ W_{\lambda, i}$ corresponds to 
$u^{\lambda}$ for all $i \in I_{\lambda}$ and $I_{\lambda}$ is the multiplicity of $u^{\lambda}$. Since 
all matrix entries  of any finite-dimensional irreducible
co-representation  of $C(SU_q(n))$ lie in the algebra generated by $\left\{u_k^l : n-m+1 \leq l \leq n, 1 \leq k \leq n \right\}$, it 
follows from  Theorem 1, \cite{PodVai-1999aa} that 
$\oplus_{\lambda \in \widehat{SU(n)}}\oplus_{i \in I_{\lambda} }W_{\lambda,i} 
\subseteq C^*\left\{u_k^l : n-m+1 \leq l \leq n, 1 \leq k \leq n \right\}$.
This proves the claim. 
\qed

We will write down  
all irreducible representations of $C(SU_q(n)/SU_q(n-m))$ 
using Theorem 2.2, \cite{NesTus-2012ab}. 
%Let $N$ be the number operator given by $N: e_{n} \mapsto n e_{n}$ 
%and $S$ be the  shift operator given by $S: e_{n} \mapsto e_{n-1}$ on $L_{2}(\bbn)$. 
For $i=1,2,\cdots,n-1$, define 
the map $\pi_{s_i}: C(SU_q(n)) \longrightarrow \cll(L_2(\bbn))$ by 
\[
\pi_{s_{i}}(u_l^k)=\begin{cases}
              \sqrt{1-q^{2N+2}}S & \mbox{ if } (k,l)=(i,i),\cr
              S^*\sqrt{1-q^{2N+2}} & \mbox{ if } (k,l)=(i+1,i+1),\cr
							-q^{N+1} & \mbox{ if } (k,l)=(i,i+1),\cr
							q^N & \mbox{ if } (k,l)=(i+1,i),\cr
							\delta_{kl} & \mbox{ otherwise }. \cr
							\end{cases}
\]
Each $\pi_{s_{i}}$ is an  irreducible representation $C(SU_{q}(n))$.
For any two representations $\varphi$ and $\psi$ of $C(SU_{q}(n))$, define $\varphi * \psi := (\varphi \otimes \psi)\circ \Delta$.
Let $W$ be the Weyl group of the Lie algebra $su_n$ and 
$\vartheta \in W$ such that  $s_{i_{1}}s_{i_{2}}...s_{i_{k}}$ is a reduced expression for $\vartheta$. 
Then $\pi_{\vartheta}= \pi_{s_{i_{1}}}*\pi_{s_{i_{2}}}*\cdots *\pi_{s_{i_{k}}}$ is an irreducible representation which is independent 
of the reduced expression. Now for $t=(t_{1},t_{2},\cdots ,t_{n-1}) \in \bbbt^{n-1}$, define the map  $\tau_t: C(SU_{q}(n)) \longrightarrow \bbc $ by
\[
\tau_{t}(u_l^k)=\begin{cases}
              t_{n-k+1}\delta_{kl} & \mbox{ if } k > 1,\cr
	      \overline{t_1t_2 \cdots t_{n-1}\delta_{kl}} & \mbox{ if } k = 1,\cr
							\end{cases}
\]
Then $\tau_{t}$ is a $C^*$-algebra homomorphism. For $t \in \bbbt^{n-1}, \vartheta \in W$,  let $\pi_{(t,\vartheta)} = \tau_{t}*\pi_{\vartheta}$. Define 
$\eta_{(t,\vartheta)}$ to be the map $\pi_{(t,\vartheta)}$ restricted to $C(SU_q(n)/SU_q(n-m))$.
Let $a=(a_1,a_2, \cdots , a_m)$ be an $m$-tuple of positive integers such that $1 \leq a_j \leq n-j+1$ for all $j \in \{1,2, \cdots , m\}$. Define 
the word  $w(a)$ of the Weyl group $W$ as follows:
\[
w(a) =s_{n-m}s_{n-m-1}\cdots s_{a_m}s_{n-m+1}s_{n-m}\cdots s_{a_{m-1}} \cdots s_{n-1}s_{n-2}\cdots s_{a_1}
 \]
 with the convention that for any $1 \leq j \leq m$, the string $s_{n-j}s_{n-j-1}\cdots s_{a_j}$ is empty  if $a_j=n-j+1$.
 Clearly $w(a)$ is in a reduced form. For $t \in \bbbt^m$, define $[t]_n$  to be $(t,1,\cdots ,1) \in \bbbt^{n-1}$. The following theorem 
 describes all irreducible representations of the $C^*$-algebra underlying the quotient  space $SU_q(n)/SU_q(n-m)$. For 
 the proof, we refer the reader to \cite{NesTus-2012ab}.
 \bthm \label{representation}
Let $\mathcal{A}_m=\left\{(a_1, \cdots , a_m): 1 \leq a_j \leq n-j+1 \mbox{ for all } j \in \{1,2, \cdots , m\}\right\}$.  
Then the set  $\{\eta_{([t]_n,w(a))}\}_{t \in \bbbt^m,a \in \mathcal{A}_m}$ is the set of all irreducible representations of the $C^*$-algebra 
 $C(SU_q(n)/SU_q(n-m))$. 
\ethm
Now we will associate some diagrams with the 
above representations
which  will be useful for our purpose.  We will use the scheme followed by 
Chakraborty \& Pal~\cite{ChaPal-2003ac} with a few additions. 
For convenience, we use labeled lines  to 
represent operators as given in the following table.
 \begin{center}
         \begin{tabular}{|c|c|c|c|}
   \hline
   Arrow type & Operator &  Arrow type  & Operator \\
   \hline   
    %&&
  %\multirow{2}{*}{}&
%\multirow{2}{*}{} & 
  %  \multirow{2}{*}{}\\
\begin{tikzpicture}[scale=.7]
\draw [-] (0,0) -- (1,0);
\end{tikzpicture}
  &
  $I$ & 
  \begin{tikzpicture}[scale=.7]
\draw [-] (0,0) -- (1,0);
\node at (.7,.3) {$t$};
\end{tikzpicture}
& $M_t$  \\ 
     \hline   
 \begin{tikzpicture}[scale=.7]
 \draw [-] (0,0) -- (1,0);
 \node at (.7,.3) {$+$};
 \end{tikzpicture}
   &
   $S^*\sqrt{I-q^{2N+2}}$ & 
\begin{tikzpicture}[scale=.7]
\draw [-] (0,0) -- (1,0);
\node at (.5,-.2) {$-$};
\end{tikzpicture}
&
  $\sqrt{I-q^{2N+2}}S$  \\
     \hline
    \begin{tikzpicture}[scale=.7]
\draw [-] (0,.7) -- (.7,0);
\end{tikzpicture}
       &
$q^{N}$ & 
\begin{tikzpicture}[scale=.7]
\draw [-] (.1,.1) -- (.8,.8);
\end{tikzpicture}
  &
$-q^{N+1}$
\\
     \hline
\end{tabular} \\ 

\end{center}

Note that for $t \in \bbbt$,  $M_t$ represents the multiplication operator   on $\bbc$ sending $1$ to $t$. For other operators, the Hilbert spaces on which they act are given at the top of the diagram. 
 Let us describe how to use a diagram to represent the irreducible representations $\pi_{s_i}$ and $\tau_t$ where $1\leq i \leq n-1$ and $t \in \bbbt^{n-1}$.
 
 \begin{tabular}{p{200pt}p{200pt}}
	% $i\neq n$:  &  $i=n$:\\
\begin{tikzpicture}[scale=1.2, shift={(0,2)}]
	 % \begin{scope}[yshift=2cm]
\draw [-] (0,3) -- (1,3);
\draw [dashed] (0,2.5) -- (1,2.5);
\draw [-] (0,2) -- (1,2);
\draw [-] (0,2) -- (1,1);
\draw [-] (0,1) -- (1,2);
\draw [-] (0,1) -- (1,1);
\draw [dashed] (0,.5) -- (1,.5);
\draw [-] (0,0) -- (1,0);
\node at (.5,2.2) {+};
\node at (1,1.7) {};
\node at (1,1.3) {};
\node at (.5,.85) {${}-{}$};
\node at (-.3,3.05){$n$};
\node at (1.3,3.05){$n$};
\node at (-.5,2.05){$i+1$};
\node at (1.5,2.05){$i+1$};
\node at (-.3,1.05){$i$};
\node at (1.3,1.05){$i$};
\node at (-.3,.05){$1$};
\node at (1.3,.05){$1$};
\node at (.5,4){$L_2(\bbn)$};
 \node at (.5,-1){\text Diagram 1: $\pi_{s_i}$};
 % \end{scope}
\end{tikzpicture}
&
\begin{tikzpicture}[scale=1.2, shift={(0,2)}]
 \draw [-] (0,3) -- (.5,3) node [above]{${}t_1$} -- (1,3);
 \draw [-] (0,2.5) -- (.5,2.5) node [above]{${}t_2$} -- (1,2.5);
 \draw [dashed] (0,2) --  (1,2);
 \draw [-] (0,1.5) -- (.5,1.5) node [above]{${}t_{m}$} -- (1,1.5);
 \draw [-] (0,1) -- (1,1);
 \draw [dashed] (0,.6) -- (1,.6);
 \draw [-] (0,0) -- (.5,0) node [above]{${}\overline{t_1}\cdots \overline{t_{m}}$} -- (1,0);
 \node at (-1,3){$n$};
 \node at (-1,2.5){$n-1$};
 \node at (-1,1.5){$n-m+1$};
 \node at (-1,1){$n-m$};
  \node at (-1,0){$1$};
 \node at (2,3){$n$};
 \node at (2,2.5){$n-1$};
  \node at (2,1.5){$n-m+1$};
 \node at (2,1){$n-m$};
  \node at (2,0){$1$};
 \node at (.5,4){$\bbc$};
  \node at (.5,-1){\text Diagram 2: $\tau_{[t]_n}$};
\end{tikzpicture}
\end{tabular}\\[1ex]		
In these two diagrams, each path from a node $k$ on the
left to a node $l$ on the right stands for an
operator  given as in the table acting on the Hilbert space given at the top of the diagrams. Now $\pi_{s_i}(u_{l}^{k})$ and $\tau_{[t]_n}(u_{l}^{k})$ are the operators represented by the path from $k$ to $l$ in 
diagram 1 and diagram 2 respectively
and are zero if there is no such path.
Thus, for example, $\pi_{s_i}(u_1^1)$
is $I$; $\pi_{s_i}(u_1^2)$ is zero whereas $\pi_{s_i}(u_{i+1}^{i})=-q^{N+1}$ if $i>1$. Similarly $\tau_{[t]_n}(u_n^n)=M_{t_1}$ and $\tau_{[t]_n}(u_n^{n-1})=0$.

 Next, let us explain how to represent $\pi \ast \rho$ by a diagram where $\pi$ and $\rho$ are two representations of  $C(SU_q(n)$ acting on the Hilbert spaces $\clh_1$ and $\clh_2$ respectively.
Simply keep the two diagrams representing $\pi$ and $\rho$
adjacent to each other. Identify, for each row, the node on the right side
of the diagram for $\pi$ with the corresponding node on the left in the diagram
for $\rho$. Now, $(\pi \ast \rho)(u_{l}^{k})$ would be
an operator on the Hilbert space $\clh_1\otimes \clh_2$ 
determined by all the paths from
the node $k$ on the left to the node $l$ on the right. It would be zero if
there is no such path and if there are more than one paths, then it would be the sum of
the operators given by each such path.
In this way, we can draw diagrams for each irreducible representation of $C(SU_q(n))$ and  $C(SU_q(n)/SU_q(n-m))$.

The following diagram  is for the representation $\eta_{([t]_n,w(a))}$ of  $C(SU_q(6)/SU_q(4))$ where $a=(3,2)$ and $t=(t_1,t_2)$.

\begin{center}
 \begin{tikzpicture}[scale=1.2]
 \draw [-] (0,6) -- (.5,6) node [above]{${}t_1$} -- (1,6)-- (2,6) -- (3,6) -- (4,6)  -- (4.5,6) node [above]{${}+$} -- (5,6) -- (6,6) -- (7,6);
 \draw [-] ((0,5) -- (.5,5) node [above]{${}t_2$} -- (1,5) -- (1.5,5) node [above]{${}+$} -- (2,5) -- (3,5) -- (4,5) -- (4.5,5) node [below]{${}-$} -- (5,5) -- (5.5,5) node [above]{${}+$} -- (6,5) -- (7,5);
 \draw [-] (0,4) -- (1,4) -- (1.5,4) node [below]{${}-$} -- (2,4) -- (2.5,4) node [above]{${}+$} -- (3,4) -- (4,4) -- (5,4) -- (5.5,4)  node [below]{${}-$} -- (6,4) -- (6.5,4) node [above]{${}+$} -- (7,4) ;
 \draw [-] (0,3) -- (1,3) -- (2,3)-- (2.5,3) node [below]{${}-$} -- (3,3) -- (3.5,3) node [above]{${}+$} -- (4,3) -- (5,3) -- (6,3) -- (6.5,3)  node [below]{${}-$} -- (7,3);
 \draw [-] (0,2) -- (1,2) -- (2,2) -- (3,2) -- (3.5,2) node [below]{${}-$} -- (4,2) -- (5,2) -- (6,2) -- (7,2);
 \draw [-] (0,1) -- (.5,1) node [above]{${}\overline{t_{1}t_2}$} -- (1,1) -- (2,1) -- (3,1) -- (4,1) -- (5,1) -- (6,1) -- (7,1);
 \draw [-] (4,6) -- (5,5) node [below]{${}\;$};
 \draw [-] (5,5) -- (6,4) node [below]{${}\;$};
 \draw [-] (6,4) -- (7,3) node [below]{${}\;$};
 \draw [-] (1,5) -- (2,4) node [below]{${}\;$};
 \draw [-] (2,4) -- (3,3) node [below]{${}\;$};
 \draw [-] (3,3) -- (4,2) node [below]{${}\;$};
 \draw [-] (1,4) -- (2,5) node [above]{${}\;$};
 \draw [-] (2,3) -- (3,4) node [above]{${}\;$};
 \draw [-] (3,2) -- (4,3) node [above]{${}\;$};
 \draw [-] (4,5) -- (5,6) node [above]{${}\;$};
 \draw [-] (5,4) -- (6,5) node [above]{${}\;$};
 \draw [-] (6,3) -- (7,4) node [above]{${}\;$};
 \node at (0,6){$\bullet$};
 \node at (1,6){$\bullet$};
 \node at (2,6){$\bullet$};
 \node at (3,6){$\bullet$};
 \node at (4,6){$\bullet$};
 \node at (5,6){$\bullet$};
 \node at (6,6){$\bullet$};
 \node at (7,6){$\bullet$};
 \node at (0,5){$\bullet$};
 \node at (1,5){$\bullet$};
 \node at (2,5){$\bullet$};
 \node at (3,5){$\bullet$};
 \node at (4,5){$\bullet$};
 \node at (5,5){$\bullet$};
 \node at (6,5){$\bullet$};
 \node at (7,5){$\bullet$};
 \node at (0,4){$\bullet$};
 \node at (1,4){$\bullet$};
 \node at (2,4){$\bullet$};
 \node at (3,4){$\bullet$};
 \node at (4,4){$\bullet$};
 \node at (5,4){$\bullet$};
 \node at (6,4){$\bullet$};
 \node at (7,4){$\bullet$};
 \node at (0,3){$\bullet$};
 \node at (1,3){$\bullet$};
 \node at (2,3){$\bullet$};
 \node at (3,3){$\bullet$};
 \node at (4,3){$\bullet$};
 \node at (5,3){$\bullet$};
 \node at (6,3){$\bullet$};
 \node at (7,3){$\bullet$};
 \node at (0,2){$\bullet$};
 \node at (1,2){$\bullet$};
 \node at (2,2){$\bullet$};
 \node at (3,2){$\bullet$};
 \node at (4,2){$\bullet$};
 \node at (5,2){$\bullet$};
 \node at (6,2){$\bullet$};
 \node at (7,2){$\bullet$};
 \node at (0,1){$\bullet$};
 \node at (1,1){$\bullet$};
 \node at (2,1){$\bullet$};
 \node at (3,1){$\bullet$};
 \node at (4,1){$\bullet$};
 \node at (5,1){$\bullet$};
 \node at (6,1){$\bullet$};
 \node at (7,1){$\bullet$};
 \node at (.7,6.7){$\bbc$};
 \node at (1.5,6.7){$L_2(\bbn)$};
 \node at (2.5,6.7){$L_2(\bbn)$};
 \node at (3.5,6.7){$L_2(\bbn)$};
 \node at (4.5,6.7){$L_2(\bbn)$};
 \node at (5.5,6.7){$L_2(\bbn)$};
 \node at (6.5,6.7){$L_2(\bbn)$};
 \node at (1,6.7){$\otimes$};
  \node at (2,6.7){$\otimes$};
  \node at (3,6.7){$\otimes$};
 \node at (4,6.7){$\otimes$};
  \node at (5,6.7){$\otimes$};
   \node at (6,6.7){$\otimes$};
   \node at (-.5,6){$6$};
   \node at (-.5,5){$5$};
   \node at (-.5,4){$4$};
   \node at (-.5,3){$3$};
   \node at (-.5,2){$2$};
   \node at (-.5,1){$1$};
   \node at (7.5,6){$6$};
   \node at (7.5,5){$5$};
   \node at (7.5,4){$4$};
   \node at (7.5,3){$3$};
   \node at (7.5,2){$2$};
   \node at (7.5,1){$1$};
 \node at (3.5,0){\text Diagram 3: $\eta_{[t]_n,w(a)}$};
 \end{tikzpicture}
 \end{center}

Let $a=(a_1,a_2,\cdots ,a_m) \in \cla_m$. For $1 \leq i \leq m$,  we denote by $w(a)_i$ the reduced word $s_{n-i}s_{n-i-1}\cdots s_{a_i}$. Therefore, 
one can write the reduced Weyl word $w(a)$ and hence the corresponding irreducible representation $\pi_{w(a)}$ of $C(SU_q(n))$  as follows. 
\[
 w(a)=w(a)_mw(a)_{m-1}\cdots w(a)_1, \qquad \qquad \pi_{w(a)}= \pi_{w(a)_m}\ast\pi_{w(a)_{m-1}}\ast \cdots \ast \pi_{w(a)_1}.
\]
The following diagram  is for the representation
$\pi_{w(a)_i}:=\pi_{s_{n-i}}\ast \pi_{s_{n-i-1}}\ast \cdots \ast \pi_{s_{a_i}}$ of $C(SU_q(n))$. 

 \begin{center}
 \begin{tikzpicture}[scale=1.2]
\draw [-] (0,10) -- (1,10) -- (2,10) -- (3,10) -- (4,10);
\draw [dashed] (0,9) -- (1,9) -- (2,9) -- (3,9) -- (4,9);
\draw [-] (0,8) -- (.5,8) node [above]{${}+$} -- (1,8) -- (2,8) -- (3,8) -- (4,8);
\draw [-] (0,7) -- (.5,7) node [below]{${}-$} -- (1,7) -- (1.5,7) node [above]{${}+$} -- (2,7) -- (3,7) -- (4,7);
\draw [-] (0,6) -- (1,6) -- (1.5,6) node [below]{${}-$} -- (2,6) -- (3,6) -- (4,6);
\draw [dashed] (0,5) -- (1,5) -- (2,5) -- (3,5) -- (4,5);
\draw [-] (0,4) -- (1,4) -- (2,4) -- (3,4) -- (3.5,4) node [above]{${}+$} -- (4,4);
\draw [-] (0,3) -- (1,3) -- (2,3) -- (3,3) -- (3.5,3) node [below]{${}-$} -- (4,3);
\draw [-] (0,2) -- (1,2) -- (2,2) -- (3,2) -- (4,2);
\draw [dashed] (0,1) -- (1,1) -- (2,1) -- (3,1) -- (4,1);
\draw [-] (0,0) -- (1,0) -- (2,0) -- (3,0) -- (4,0);
\node at (0,10){$\bullet$};
 \node at (1,10){$\bullet$};
 \node at (2,10){$\bullet$};
 \node at (3,10){$\bullet$};
 \node at (4,10){$\bullet$};
 \node at (0,8){$\bullet$};
 \node at (1,8){$\bullet$};
 \node at (2,8){$\bullet$};
 \node at (3,8){$\bullet$};
 \node at (4,8){$\bullet$};
 \node at (0,7){$\bullet$};
 \node at (1,7){$\bullet$};
 \node at (2,7){$\bullet$};
 \node at (3,7){$\bullet$};
 \node at (4,7){$\bullet$};
 \node at (0,6){$\bullet$};
 \node at (1,6){$\bullet$};
 \node at (2,6){$\bullet$};
 \node at (3,6){$\bullet$};
 \node at (4,6){$\bullet$};
 \node at (0,5){$\bullet$};
 \node at (1,5){$\bullet$};
 \node at (2,5){$\bullet$};
 \node at (3,5){$\bullet$};
 \node at (4,5){$\bullet$};
 \node at (0,4){$\bullet$};
 \node at (1,4){$\bullet$};
 \node at (2,4){$\bullet$};
 \node at (3,4){$\bullet$};
 \node at (4,4){$\bullet$};
 \node at (0,3){$\bullet$};
 \node at (1,3){$\bullet$};
 \node at (2,3){$\bullet$};
 \node at (3,3){$\bullet$};
 \node at (4,3){$\bullet$};
 \node at (0,2){$\bullet$};
 \node at (1,2){$\bullet$};
 \node at (2,2){$\bullet$};
 \node at (3,2){$\bullet$};
 \node at (4,2){$\bullet$};
 \node at (0,0){$\bullet$};
 \node at (1,0){$\bullet$};
 \node at (2,0){$\bullet$};
 \node at (3,0){$\bullet$};
 \node at (4,0){$\bullet$};
 \draw [-] (0,8) -- (1,7) node [below]{${}\;$};
  \draw [-] (0,7) -- (1,8) node [above]{${}\;$};
  \draw [-] (1,7) -- (2,6) node [below]{${}\;$};
  \draw [-] (1,6) -- (2,7) node [above]{${}\;$};
  \draw [-] (3,4) -- (4,3) node [below]{${}\;$};
  \draw [-] (3,3) -- (4,4) node [above]{${}\;$};
  \node at (-1,10){$n$};
  \node at (-1,8){$n-i+1$};
  \node at (-1,7){$n-i$};
  \node at (-1,6){$n-i$};
  \node at (-1,4){$a_i+1$};
  \node at (-1,3){$a_i$};
  \node at (-1,0){$1$};
  \node at (5,10){$n$};
  \node at (5,8){$n-i+1$};
  \node at (5,7){$n-i$};
  \node at (5,6){$n-i$};
  \node at (5,4){$a_i+1$};
  \node at (5,3){$a_i$};
  \node at (5,0){$1$};
  \node at (.5,10.7){$L_2(\bbn)$};
  \node at (1.5,10.7){$L_2(\bbn)$};
 \node at (2.5,10.7){$\cdots$};
 \node at (3.5,10.7){$L_2(\bbn)$};
 \node at (1,10.7){$\otimes$};
  \node at (2,10.7){$\otimes$};
  \node at (3,10.7){$\otimes$};
   \node at (2,-1){\text Diagram 4: $\pi(w(a)_i)$};
  \end{tikzpicture}
 \end{center}

\newsection{Quantum Stiefel manifolds}

 In this section, our aim is to  establish $C(SU_q(n)/SU_q(n-m))$  as a universal $C^*$-algebra given by 
 finite sets of generators 
 and relations.  We start with a definition.
 
\bdfn 
  We define $C^{*}$-algebra $C^{n,m}$ as the universal  
 $C^{*}$-algebra generated by the elements $\{w_j^i: n-m+1 \leq i \leq n , 1\leq j \leq n\}$ satisfying the following relations:
\begin{IEEEeqnarray}{rCll} 
 w_{k}^iw_{k}^j &=&  q w_{k}^jw_{k}^i  &\mbox{for }  i<j,  \label{c1}\\
 w_{k}^iw_{l}^i &=& qw_{l}^iw_{k}^i  & \mbox{for }  k<l,  \label{c2}\\
 w_{k}^iw_{l}^j &=& w_{l}^jw_{k}^i  & \mbox{for }  i<j, k>l, \label{c3}\\
 w_{k}^iw_{l}^j &= &w_{l}^jw_{k}^i+(q^{-1}-q)w_{l}^iw_{k}^j  &\mbox{for }  i>j, k>l,\label{c4}\\ 
 (w_{k}^i)^*w_{l}^j&=&w_{l}^j(w_{k}^i)^* &  \mbox{for } i\neq j, k\neq l, \label{c5}\\
 (w_{k}^i)^*w_{l}^i+(1-q^2)\sum_{j>i}(w_{k}^j)^*w_{l}^j&=&qw_{l}^i(w_{k}^i)^* &\mbox{for } k \neq l, \label{c6}\\
 w_{k}^i(w_{k}^j)^*+(1-q^2)\sum_{l<k}w_{l}^i(w_{l}^j)^*&=&q(w_{k}^i)^*w_{k}^j &\mbox{for } i \neq j, \label{c7}\\
 (w_{k}^i)^*w_{k}^i+(1-q^2)\sum_{j>i}(w_{k}^j)^*w_{k}^j &=&w_{k}^i(w_{k}^i)^*+(1-q^2)\sum_{l<k}w_{l}^i(w_{l}^i)^* \label{c8}\\
 \sum_{k=1}^nw_{k}^i(w_{k}^j)^*&=&\delta_{ij}. &\label{c9}
\end{IEEEeqnarray}
\edfn
\brmrk
These relations will be called commutation relations. 
It can be easily obtained from the relations given in \cite{PodVai-1999aa} or using $R$-matrix for type $A_{n-1}$ 
(see \cite{KliSch-1997aa}). Hence the generators $\{u_j^i: n-m+1 \leq i \leq n , 1\leq j \leq n\}$ of 
 $C(SU_q(n)/SU_q(n-m))$ satisfy these relations. Using this fact and 
relation
(\ref{c9}), one can  prove that the $C^*$-algebra $C^{n,m}$ is well defined.
\ermrk Since the generators $\{u_j^i: n-m+1 \leq i \leq n , 1\leq j \leq n\}$  of $C(SU_q(n)/SU_q(n-m))$ satisfy the same relations 
satisfied by the generators $\{w_j^i: n-m+1 \leq i \leq n , 1\leq j \leq n\}$ of $C^{n,m}$, it follows from the universal 
property of $C^{n,m}$ that there exists a unique surjective homomorphism 
\[
 \Psi_m:C^{n,m} \rightarrow C(SU_q(n)/SU_q(n-m))
\]
 sending 
$w_j^i$ to $u_j^i$ for all $1\leq j \leq n$ and $n-m+1 \leq i \leq n $. We will show that for all $1\leq m \leq n$, the map $\Psi_m$ is an isomorphism. For that, we need   
  to analyse an irreducible representation of the $C(SU_q(n)/SU_q(n-m))$ more closely.

Let  $\pi$ be an irreducible representation of $C(SU_q(n)/SU_q(n-m))$ on the Hilbert space $\clh$. Then for 
some $t=(t_1,\cdots ,t_m) \in \bbbt^m$ and $a=(a_1,\cdots ,a_m) \in \cla_m$, the representation  $\pi$ is 
equivalent to the representation $\eta_{([t]_n,w(a))}$.   Set
\[
 \clh_0^1=\clh, \quad F_1=\{1,2,\cdots n\}, \quad \ell_1=n-a_1+1, \quad c_1=a_1. 
\]
Let $T_k^1=\pi(u_{n-k+1}^n) \mbox{ for } 1 \leq k \leq n$ acting on the Hilbert space $\clh_0^1$. 
Using  commutation relations, one can check  that  
the operators $(T_1^1, \cdots T_{n}^1)$ acting on $\clh_0^1$ satisfy $ODQS(q)$ relations with rank $\ell_1$.   Suppose that  $\clh_0^{j}, F_j, \ell_j, c_j$ and the operators $\{T_k^j:\mbox{ for }
  1\leq k \leq n-j+1\}$ on $\clh_0^{j}$ is defined for all $j=1,2,\cdots i-1$. Further assume that for  $j=1,2,\cdots i-1$,
  the ordered tuple $(T_1^j, \cdots T_{n-j+1}^j)$ satisfy $ODQS(q)$ relations with rank $\ell_j$. Define
\begin{IEEEeqnarray}{rCl} \label{c}
\clh_0^{i}&=& \mbox{ eigenspace of } T_{\ell_{i-1}}^{i-1} \mbox{ corresponding to eigenvalue } 1,\nonumber  \\ 
\ell_i&=&n+2-i-a_i, \nonumber \\
F_i&=&\{1,2,\cdots,n\}-\{c_1,c_2,\cdots , c_{i-1}\}, \nonumber \\
g_i(f)&=&\#\{j:c_j<f\} \quad \qquad \quad  \mbox{ for } f \in F_i,\nonumber \\
c_i&=&f_{\ell_i}^i,
\end{IEEEeqnarray}

\begin{equation} \label{operators}
 T_k^i= \frac{(-1)^{g_{i}(f_k^{i})}}{q^{g_{i}(f_k^{i})}}\pi(u_{f_k^{i}}^{n-i+1})_{|\clh_0^{i}} \quad \mbox{ for }
  1\leq k \leq n-i+1,
\end{equation}
where $f_1^{i}>f_2^{i}>\cdots > f_{n-i+1}^{i}$ are all elements of $F_i$. In this way, we define the operators $T_1^i, T_2^i, \cdots, T_{n-i+1}^i$ for all $i \in \{1,2,\cdots, m\}$. One thing we need to show that at each stage, 
the tuple $(T_1^i,T_2^i,\cdots,T_{n-i+1}^i)$  satisfies ODQS(q) relations with rank $\ell_i$. We will show this by identifying these operators using the diagram associated to the  irreducible representation $\eta_{([t]_n,w(a))}$.
 First note the following facts from diagram 4.  
\begin{enumerate}
 \item 
 The eigenspace of the operator $\pi(w(a)_i)(u_{a_i}^{n-i+1}(u_{a_i}^{n-i+1})^*)$ corresponding to eigenvalue $1$ is one dimensional subspace $K_0^i$ spanned by the vector $e_0\otimes e_0 \otimes \cdots \otimes e_0$. 
 \item 
 Let $l < n-i+1$. Then on $K_0^i$, starting from $l$, one can go to either $l+1$ or $l$ depending on whether $l\geq a_i$ or $l<a_i$. Moreover, in this case we have 
 \begin{IEEEeqnarray}{rCll}
  \pi(w(a)_i)(u_l^{l+1})_{|K_0^i} &=& qI \quad &\mbox { if } l\geq a_i \nonumber \\
  \pi(w(a)_i)(u_l^{l})_{|K_0^i} &=& I \quad &\mbox { if } l< a_i \nonumber 
 \end{IEEEeqnarray}
\end{enumerate}
Using these observations and the diagram associated with the representation $\eta_{([t]_n,w(a))}$, one can show that on the Hilbert space $\clh_0^i$, $T_{\ell_i}^i=t_{n-i+1}I$ in  case $\ell_i=1$ and in  case  $\ell_i>1$, we have 
 \[
  T_l^i= \begin{cases}
          t_{n-i+1}\underbrace{1 \otimes 1 \otimes \cdots \otimes 1}_{\sum_{j=i+1}^m(\ell_j-1) \mbox{ copies}}\otimes \sqrt{1-q^{2N}}S^* \otimes 1 \otimes 1 \otimes \cdots \otimes 1 
          & \mbox{ if }  l=1; \cr
          t_{n-i+1}\underbrace{1 \otimes 1 \otimes \cdots \otimes 1}_{\sum_{j=i+1}^m(\ell_j-1) \mbox{ copies}}\otimes \underbrace{q^N \otimes  \cdots \otimes q^N}_{l-1 \mbox{ copies} } 
          \otimes  \sqrt{1-q^{2N}}S^* \otimes 1 \otimes \cdots \otimes 1 
          & \mbox{ if }  1< l <\ell_i; \cr
          t_{n-i+1}\underbrace{1 \otimes 1 \otimes \cdots \otimes 1}_{\sum_{j=i+1}^m(\ell_j-1) \mbox{ copies}}\otimes \underbrace{q^N \otimes  \cdots \otimes q^N}_{l-1 \mbox{ copies} } 
          \otimes 1 \otimes 1 \otimes \cdots \otimes 1 
          & \mbox{ if }  l =\ell_i; \cr
          0 & \mbox{ if }  l >\ell_i ; \cr
         \end{cases}
 \]

 Hence  for $1\leq i\leq m$, the tuple $(T_1^i,T_2^i,\cdots,T_{n-i+1}^i)$  satisfies ODQS(q) relations with rank $\ell_i$. Moreover, 
 \begin{enumerate}
  \item 
   The eigenspace $\clh_0^{m+1}$ of the operator $T_{\ell_m}^m(T_{\ell_m}^m)^*$ corresponding 
   to eigenvalue $1$
 is one dimensional. It follows from part (10) of the proposition \ref{ppsn-spectrum} that the tuple $(T_1^m,T_2^m,\cdots,T_{n+1-m}^m)$ is irreducible.
\item
For any nonzero vector $h \in \clh_0^{m+1}$,  $T_{\ell_i}^ih=t_ih$ 
for $1 \leq i \leq m$.
\item 
Let $1\leq j \leq m$. Then on the Hilbert space $\clh_0^j$,  
\begin{equation}\label{kerc_j=0}
\ker(T_{\ell_j}^j(T_{\ell_j}^j)^*)=\ker(T_{\ell_j}^j)=\{0\} 
\end{equation}

\item  Let $k<j\leq m$. Then on the Hilbert space $\clh_0^{j}$, one has
\begin{equation} \label{c_k=0}
  \pi(u_{c_k}^{n-j+1})=0 
\end{equation} 
 \item 
  Let  $j \leq m$ and $k<c_j$. Then on the Hilbert space $\clh_0^{j}$, one has
  \begin{equation} \label{<c_k}
   \pi(u_{k}^{n-j+1})=0
  \end{equation}
 
\end{enumerate} 
 \brmrk \label{define T} For any $1\leq i \leq m$, the function $g_i$ and the set $F_i$ depends on the set $\{a_1,a_2,\cdots a_{i-1}\}$ only. Therefore, the expression 
 of the operators $T_1^i,\cdots T_{n-i+1}^i$ depends on $\{a_1,a_2,\cdots a_{i-1}\}$ only and hence  we can   
  define the operators $T_1^i,\cdots T_{n-i+1}^i$ for the representation $\pi$ of $C(SU_q(n)/SU_q(n-m)$ if $\pi$
 can be decomposed into those irreducible representations $\eta_{([t]_n,w(a))}$'s which have same value of $\{a_1,a_2,\cdots a_{i-1}\}$. In addition to this, one can define   $\clh_0^j$'s and  $c_j$'s for $j < i$ in the same manner.
 \ermrk

%$D_i^{m,k}=\mbox{min}\{n+2-m-k,a_{m-1},\cdots,a_{i}\}$. 
%Set $T_1^i=t_{n-i+1}$ if $a_i=1$. For $ i \in \{1,2,\cdots, m\}$ such that $a_i >1$, 
%We will define inductively
%the operators $T_1^i, T_2^i, \cdots T_{n+1-i}^i$  for $1 \leq i \leq m$ as follows: Set $T_k^1=\pi(u_k^n)$ for $i \in \{1,2,\cdots,n\}$.
%Having defined $T_1^i, T_2^i, \cdots T_{n+1-i}^i$  for $1 \leq i \leq l-1$, we define  for $ k \in \{1,2,\cdots, n+1-l\}$,  
%\begin{IEEEeqnarray}{rCl}\label{operator}
%T_k^l&=&\pi(u_{D_{1}^{l,k}}^{n-l+1})\prod_{i=1}^{l-1} \frac{1}{1-q^2} T_{n+1-i-D_i^{l,k}}^i
%\end{IEEEeqnarray}
%with the convention that if $D_i^{l,k} > D_{i-1}^{l,k}$ then instead of the operator $\frac{1}{1-q^2} T_{n+1-i-D_i^{l,k}}^i$, 
%we take the operator $1$. 

To show that $C(SU_q(n)/SU_q(n-m))$ is isomorphic to the $C^*$-algebra $C^{n,m}$, we use induction on $m$.
Fix $n$. We first take up the case of  $m=1$. Although the following proposition 
has been already proved (see \cite{VakSou-1990ab}), we prove it for the sake of completeness. 
\bppsn \label{m=1}
The map  $\Psi_1:C^{n,1} \rightarrow C(SU_q(n)/SU_q(n-1))$ sending $w_j^n$ to $u_j^n$ for all $1\leq j \leq n$ is an isomorphism.
\eppsn
\prf
It is enough to show that the map $\Psi_1$ is injective. Take an irreducible representation $\pi$ 
of $C^{n,1}$. Using  commutation relations, one can easily see that 
the tuple $(\pi(w_1^n),\pi(w_2^n),\cdots,\pi(w_n^n))$ satisfy ODQS(q) relations. Let $l_0$ and $t_0$ be the rank and the angle respectively
of this ordered tuple. Define $a=n-l_0+1$.  Consider the ordered tuple 
$(\eta_{(\tilde{t_0},w(a))}(u_1^n), \eta_{(\tilde{t_0},w(a))}(u_2^n), \cdots,\eta_{\tilde{(t_0},w(a))}(u_n^n))$.
It is easy to check that rank and angle of this ordered tuple is $l_0$ and $t_0$ respectively. Hence by proposition \ref{isomordered}, these 
ordered tuples are isomorphic which further implies that the representations $\pi$ and $\eta_{(\tilde{t_0},w(a))}\circ \Psi_1$  of $C^{n,1}$ are equivalent. 
This shows that any irreducible representation of $C^{n,1}$ factors through the map $\Psi_1$. Hence image of an element in $\ker(\Psi_1)$ 
under any irreducible representation of $C^{n,1}$ is zero which implies that that element is zero. This proves the claim.
\qed \\
Suppose that for $m=1,2,\cdots k-1$, the map $\Psi_m$ is an isomorphism between   $C^{n,m}$ and   $C(SU_q(n)/SU_q(n-m))$. Let $\pi$ 
be an irreducible representation of $C^{n,k}$ acting on a Hilbert space $\clh$. Let $\eta$ be the representation $\pi$ restricted to the $C^*$-algebra  generated by $\{w_j^i: n-k+2\leq i \leq n, 1\leq j \leq n\}$ that is isomorphic to 
$C^{n,k-1}$. By induction 
hypothesis and by Theorem \ref{representation}, we have
\begin{IEEEeqnarray}{rCl} \label{decomposition}
 \eta \equiv \oplus_{(t,a) \in  \bbbt^{k-1}  \times \mathcal{A}_{k-1}} \oplus_{j\in J_{(t,a)}} \eta_{([t]_n,w(a))}, \quad
 \clh \equiv \oplus_{(t,a) \in  \bbbt^{k-1} \times \mathcal{A}_{k-1}} \oplus_{j\in J_{(t,a)}}\clh_{(t,a)}^j.
\end{IEEEeqnarray}
Here $J_{(t,a)}$ represents an index set of the multiplicity of the irreducible representation $\eta_{([t]_n,w(a))}$ in the representation 
$\eta$ and for each $j \in J_{(t,a)}$, 
$\clh_{(t,a)}^j$ denotes the Hilbert space for the irreducible representation $\eta_{([t]_n,w(a))}$. We will show that $J_{(t,a)}$ 
is empty set for 
all but one value of $a$. More precisely, 
\blmma  \label{onea}
Suppose that  the map $\Psi_m:C^{n,m} \rightarrow C(SU_q(n)/SU_q(n-m))$ sending $w_j^i$ to $u_j^i$  is an isomorphism  for all $m=1,2,\cdots k-1$. Then for any irreducible
representation $\pi$ of $C^{n,k}$, the index set $J_{(t,a)}$ defined as above 
is empty set for 
all but one value of $a \in \cla_{k-1}$.  
\elmma 
\prf
Define the operators $T_i^1:=\pi(w_i^n)$ for $1\leq i \leq n$.
From commutation relations, it follows that  the ordered tuple $(T_1^1,T_2^1,\cdots,T_n^1)$ satisfies ODQS(q) relations. Let its rank be $\ell_1$ and let $a_1^0:=n-\ell_1+1$. This 
implies that for $a \in \cla_{k-1}$ such that  $a_1< a_1^0$,  $J_{(t,a)}$ is empty set. Also from the equation (\ref{decomposition}), we have
\[
 \ker(T_{\ell_1}^1)=\oplus_{\{(t,a) \in  \bbbt^{k-1} \times \mathcal{A}_{k-1}: a_1> a_1^0\}} \oplus_{j\in J_{(t,a)}}\clh_{(t,a)}^j.
\]
Clearly $\ker(T_{\ell_1}^1)$ is invariant under the action of $\{\pi(w_j^i),\pi((w_j^i)^*): 1 \leq j \leq n, n+2-k \leq i \leq n\}$. 
Using commutations relations,
one can  show that $\ker(T_{\ell_1}^1)$ is an invariant subspace of $\pi$. Since $T_{\ell_1}^1$ is nonzero operator, we have $\ker(T_{\ell_1}^1)=\{0\}$.
This 
implies  that if $a_1>a_1^0$ then $J_{(t,a)}$ is empty set.

Assume that for $a \in \cla_{k-1}$ such that 
$(a_1,a_2,\cdots,a_{i-1}) \neq (a_1^0,a_2^0, \cdots a_{i-1}^0)$,  $J_{(t,a)}$ is empty set.  
Define the operators $T_1^i,T_2^i,\cdots T_{n-i+1}$ by equation (\ref{operators}) acting on the the Hilbert space
$\clh_0^i$ (see remark \ref{define T}). Clearly the tuple $(T_1^i,T_2^i,\cdots T_{n-i+1})$ obeys ODQS(q) relations as it is a direct sum of  some tuples of operators satisfying ODQS(q) relations.  
Let its rank be $\ell_i$ and  $a_i^0:=n+2-i-\ell_i$. This implies that for $a \in \cla_{k-1}$ such that  $a_i< a_i^0$,  $J_{(t,a)}$ is empty set. Let $H_i:=\ker(T_{\ell_i}^i)$. Then using equations (\ref{decomposition}), (\ref{kerc_j=0}) and 
part (9) of the proposition \ref{ppsn-spectrum}, we get 
\begin{IEEEeqnarray}{lCl}
 \oplus_{\{(t,a) \in  \bbbt^{k-1} \times \mathcal{A}_{k-1}: a_i> a_i^0\}}\oplus_{j\in J_{(t,a)}}\clh_{(t,a)}^j 
  &=&\oplus_{\alpha_r^s \in \bbn} (T_1^1)^{\alpha_1^1}(T_2^1)^{\alpha_2^1}\cdots (T_{\ell_1-1}^1)^{\alpha_{\ell_1-1}^1} \cdots  \nonumber \\
  && \cdots (T_1^{i-1})^{\alpha_1^{i-1}}(T_2^{i-1})^{\alpha_2^{i-1}}\cdots (T_{\ell_{i-1}-1}^{i-1})^{\alpha_{\ell_{i-1}-1}^{i-1}}
  H_i. \nonumber
\end{IEEEeqnarray}
Let $H_{>a_i^0}:=\oplus_{\{(t,a) \in  \bbbt^{k-1} \times \mathcal{A}_{k-1}: a_i> a_i^0\}}\oplus_{j\in J_{(t,a)}}\clh_{(t,a)}^j$. We will show that $H_{>a_i^0}$ is an invariant subspace of $\pi$.
From the equation (\ref{decomposition}), it follows that  $H_{>a_i^0}$ is invariant under the action of $\{\pi(w_j^l),\pi(w_j^l)^*: 1 \leq j \leq n, n+2-k \leq l \leq n\}$. Using commutation relations  and equation (\ref{<c_k}),
one can show that $H_i$ is invariant under $\{\pi(w_j^{n-k+1}),\pi(w_j^{n-k+1})^*: 1 \leq j \leq n\}$. By applying  relations (\ref{c1}), (\ref{c3}), 
(\ref{c4}), (\ref{c5}), and (\ref{c7}), it follows that  $H_{>a_i^0}$  an invariant $\pi(w_n^{n-k+1})$ and $\pi(w_n^{n-k+1})^*$. 
Assume that $H_{>a_i^0}$ is invariant under the action of $\{\pi(w_j^{n-k+1}),\pi(w_j^{n-k+1})^*: r+1 \leq j \leq n\}$. 
Now using relations (\ref{c1}), (\ref{c3}), 
(\ref{c4}), (\ref{c5}), and (\ref{c7}), it follows that  $H_{>a_i^0}$ is invariant under the action of $\pi(w_r^{n-k+1})$ and $\pi(w_r^{n-k+1})^*$. Hence by backward induction, one can show 
that $H_{>a_i^0}$ is an invariant subspace of $\pi$. Since $T_{\ell_i}^i$ is nonzero operator, $H_{>a_i^0}\neq \clh$ and hence $H_{>a_i^0}=\{0\}$. This proves the claim.
\qed
\brmrk \label{ker}
It follows from the above proof that  $\ker(T_{\ell_j}^j(T_{\ell_j}^j)^*)=\ker(T_{\ell_j}^j)=\{0\}$ for all $1\leq j \leq k-1$. 
\ermrk
Define $a_0:=(a_1^0,a_2^0,\cdots,a_{k-1}^0)$. For $1\leq j \leq k-1$, let  $c_j^0$ be as defined in the equation (\ref{c}) for the representation $\eta_{([t]_n,w(a_0))}$ of $C(SU_q(n)/SU_q(n-k+1))$ for any $t \in \bbbt^{k-1}$.
Let $\clh_0^{k}$  be the eigenspace of 
$T_{\ell_{k-1}}^{k-1}(T_{\ell_{k-1}}^{k-1})^*$ corresponding to the eigenvalue $1$. Define 
\begin{IEEEeqnarray}{rCll} \label{lastoperator}
F_k&=&\{1,2,\cdots,n\}-\{c_1^0,c_2^0,\cdots , c_{k-1}^0\}, \qquad &
g_k(f)=\#\{j:c_j^0<f\} \quad \mbox{ for } f \in F_k,  \nonumber\\
 T_i^k&=& \frac{(-1)^{g_{k}(f_i^{k})}}{q^{g_{k}(f_i^{k})}}\pi(w_{f_i^{k}}^{n-k+1})_{|\clh_0^{k}} & \mbox{ for }
  1\leq i \leq n-k+1.  \nonumber
\end{IEEEeqnarray}

\bppsn \label{odqs} Let  the map $\Psi_m:C^{n,m} \rightarrow C(SU_q(n)/SU_q(n-m))$ sending $w_j^i$ to $u_j^i$  be an isomorphism for all $m=1,2,\cdots k-1$. Then  
the operators $T_1^k,T_2^k,\cdots,T_{n-k+1}^k$    defined  above are well defined.  Moreover, the tuple $(T_1^k,T_2^k,\cdots,T_{n-k+1}^k)$  satisfies ODQS(q) relations and 
 is irreducible. 
\eppsn
\prf It follows from Lemma \ref{onea} that  for each $1\leq i <k$, the tuple $(T_1^i,T_2^i,\cdots,T_{n-i+1}^i)$  satisfies ODQS(q) relations and of rank $\ell_i:=n+2-i-a_i^0$. Therefore, 
to show $T_1^k,T_2^k,\cdots,T_{n-k+1}^k$ are well 
defined operators on $\clh_0^{k}$, we need to show that for $1\leq i \leq n-k+1$ and   $h \in \clh_0^{k}$,  $T_i^k h \in \clh_0^{k}$ or equivalently $T_{\ell_{j}}^{j}T_i^k h=T_i^k h$ for all $j< k$. For all $i$ such that 
$f_i^k>f_{\ell_{j}}^{j}$, it follows directly  
from the relation (\ref{c3}).  For $i$ such that 
$f_i^k<f_{\ell_{j}}^{j}$, we have 
\begin{IEEEeqnarray}{lCl}
 \pi(w_{f_{\ell_{j}}^{j}}^{n-j+1})\pi(w_{f_i^k}^{n-k+1})h \nonumber \\
 =\pi(w_{f_i^k}^{n-k+1})\pi(w_{f_{\ell_{j}}^{j}}^{n-j+1})h+(q^{-1}-q)\pi(w_{f_{\ell_{j}}^{j}}^{n-k+1})\pi(w_{f_i^k}^{n-j+1})h  \qquad 
  (\mbox{from relation } (\ref{c4})) \nonumber \\
 =\pi(w_{f_i^k}^{n-k+1})\pi(w_{f_{\ell_{j}}^{j}}^{n-k})h  \qquad 
  (\mbox{from equation } (\ref{<c_k})) \nonumber \\
 =\pi(u_{f_i^k}^{n-k+1})h \qquad (\mbox {since } \clh_k \subset \clh_{j+1})\nonumber 
\end{IEEEeqnarray}
This shows that  $T_{\ell_{j}}^{j}T_i^k h=T_i^k h$ for $h \in \clh_0^{k}$ and $j<k$.  
Now 
we shall show that the tuple $(T_1^k,T_2^k,\cdots,T_{n-k+1}^k)$  satisfies ODQS(q) relations. From  relation 
(\ref{c1}), it follows that $T_i^kT_j^k=qT_j^kT_i^k$ for $i<j$.   
It follows from part (11) of the proposition 
\ref{ppsn-spectrum} that for  
$h \in \clh_0^k \subset \clh_0^{j+1}$ and $i \neq l$,  one has $(T_i^j)^*T_l^jh=qT_l^j(T_i^j)^*h=0$. Using this fact and the relation (\ref{c6}), we get $(T_i^k)^*T_l^k=qT_l^k(T_i^k)^*$ for $i \neq l$. 
Observe that for  $h \in \clh_0^{k}$, we have
\begin{IEEEeqnarray}{rCll}
 \frac{1}{q}\pi(u_{c_i}^{n-k+1})h&=&\frac{1}{q}\pi(u_{c_i}^{n-k+1})T_{\ell_i}^{i}h & \quad (\mbox{since } \clh_0^{k}\subset  \clh_0^{i}  )\nonumber \\ 
 &=&T_{\ell_i}^{i}\pi(u_{c_i}^{n-k+1})h & \quad (\mbox{from equation } (\ref{c1}) \mbox{ and } c_i=f_{\ell_i}^i) \nonumber \\
 &=& 0 &   \quad (\mbox{from  part (6) of the proposition } \ref{ppsn-spectrum}) \nonumber
\end{IEEEeqnarray}
Hence on  $\clh_0^{k}$, we have $\pi(u_{c_i}^{n-k+1})=0$. Let us assume that 
$c_{i_1}<c_{i_2}<\cdots c_{i_{k-1}}$ where $\{i_1,i_2,\cdots ,i_{k-1}\}=\{1,2,\cdots,k-1\}$. Then on the Hilbert space $\clh_0^k$, we have 
\begin{IEEEeqnarray}{lCl}
 \pi(w_{c_{i_{r}}}^{n-k+1})\pi(w_{c_{i_{r}}}^{n-k+1})^*\nonumber \\
 =\pi(w_{c_{i_{r}}}^{n-k+1})^*\pi(w_{c_{i_{r}}}^{n-k+1})+(1-q^2)\sum_{j>n-k+1}\pi(w_{c_{i_{r}}}^{j})^*\pi(w_{c_{i_{r}}}^{j})\nonumber \\
-(1-q^2)\sum_{l<{c_{i_{r}}}}\pi(w_l^{n-k+1})\pi(w_l^{n-k+1})^* \nonumber  \\
=(1-q^2)\sum_{j\in \{n-i_{r-1}+1,\cdots n-i_1+1\}}\pi(w_{c_{i_{r}}}^{j})^*\pi(w_{c_{i_{r}}}^{j})+(1-q^2)\pi(w_{c_{i_{r}}}^{n-i_{r}+1})^*\pi(w_{c_{i_{r}}}^{n-i_{r}+1})\nonumber \\
+(1-q^2)\sum_{j\in \{n-i_{k-1}+1,\cdots n-i_{r+1}+1\}}\pi(w_{c_{i_{r}}}^{j})^*\pi(w_{c_{i_{r}}}^{j})-(1-q^2)\sum_{l<{c_{i_{r}}}}\pi(w_l^{n-k+1})\pi(w_l^{n-k+1})^*\nonumber  \\
=(1-q^2)^2(1+q^2+\cdots q^{2k-4})+(1-q^2)q^{2k-2}-(1-q^2)\sum_{l<{c_{i_{r}}}}\pi(w_l^{n-k+1})\pi(w_l^{n-k+1})^* \nonumber  \\
\qquad  \qquad \qquad \quad \quad (\mbox{by equations } (\ref{operators}), (\ref{<c_k}) \mbox{ and  part (11) of proposition } \ref{ppsn-spectrum}) \nonumber  \\
=(1-q^2)-(1-q^2)\sum_{l<{c_{i_{r}}}}\pi(w_l^{n-k+1})\pi(w_l^{n-k+1})^*  \label{condition3}
\end{IEEEeqnarray}
A little algebra  along with the relations (\ref{c8}), (\ref{c9}) and equation (\ref{condition3}) will show that 
\begin{IEEEeqnarray}{rCl}
 (T_i^k)^*T_i^k=T_i^k(T_i^k)^*+(1-q^2)\sum_{j>i}T_j^k(T_j^k)^*, \qquad \sum_{i=1}^{n-k+1}T_i^k(T_i^k)^*=1  \nonumber 
\end{IEEEeqnarray}
Let $\ell_k$ be the rank of $(T_1^k,T_2^k,\cdots,T_{n-k+1}^k)$. Let $\clh_{k+1}^0$  be  the eigenspace of the operator $T_{\ell_k}^k(T_{\ell_k}^k)^*$ corresponding 
   to eigenvalue $1$. To show that the  tuple $(T_1^k,T_2^k,\cdots,T_{n-k+1}^k)$ is irreducible, take a nonzero vector $h_0 \in \clh_{k+1}^0$. Define 
\[K_{h_0}:=\oplus_{\alpha_j^i \in \bbn}(T_1^1)^{\alpha_1^1}(T_2^1)^{\alpha_2^1}\cdots (T_{\ell_1-1}^1)^{\alpha_{\ell_1-1}^1}\cdots (T_1^k)^{\alpha_1^k}(T_2^k)^{\alpha_2^k}\cdots 
  (T_{\ell_k-1}^k)^{\alpha_{\ell_k-1}^k}h_0
  \]
  From the equation (\ref{decomposition}), it follows that the subspace $K_{h_0}$ is invariant under the action of $\{\pi(w_j^l),\pi(w_j^l)^*: 1 \leq j \leq n, n+2-k \leq l \leq n\}$. 
  By applying commutation relations and the fact that $\pi(w_{c_i}^{n-k+1})=0$ on $\clh_0^k$, one can prove that $K_{h_0}$ is invariant under the action of $\pi(w_n^{n-k+1})$. Using backward induction and commutation relations, 
  it follows that $K_{h_0}$ is invariant under the action of $\{\pi(w_i^{n-k+1}): i=1,2,\cdots ,n\}$. In the same way, one can show that  for any $h \in \clh_0^{k+1}$, $K_h$ defined in the same way as $K_{h_0}$ is  
  invariant under the action of $\{\pi(w_i^{n-k+1}): i=1,2,\cdots ,n\}$. This implies that $\{\pi(w_i^{n-k+1}): i=1,2,\cdots ,n\}$ keep  $K_{h_0}$ as well as its complement invariant. As a consequence, the subspace $K_{h_0}$ 
  is invariant under the action of $\{\pi(w_i^{n-k+1}),\pi(w_i^{n-k+1})^*: i=1,2,\cdots ,n\}$ and hence invariant under $\pi$. Since $\pi$ is irreducible, it follows that $\clh_0^{k+1}$ is one dimensional. By part (10) of 
  the proposition \ref{ppsn-spectrum}, we get the claim.
\qed

In Lemma \ref{onea}, we  proved that any irreducible representation $\pi$ of $C^{n,k}$ when restricted to $C^{n,k-1}$ decomposes into those  irreducible 
representations of $C(SU_q(n)/SU_q(n-k+1))$  that have same value of $a \in \cla_{k-1}$. But it does not rule out the possibility that 
these irreducible representations may have different values of $t \in \bbbt^{k-1}$. In the following lemma, we will show that in such decomposition, 
only one irreducible representation of $C(SU_q(n)/SU_q(n-k+1))$ occurs with certain (infinite) multiplicity.
 
\blmma \label{one(a,t)}
Suppose that  the map $\Psi_m:C^{n,m} \rightarrow C(SU_q(n)/SU_q(n-m))$ sending $w_j^i$ to $u_j^i$  is an isomorphism  for all $m=1,2,\cdots k-1$. Then for any irreducible
representation $\pi$ of $C^{n,k}$, the index set $J_{(t,a)}$ defined as above 
is empty set for 
all but one value of $(t,a) \in \bbbt^{k-1} \times \cla_{k-1}$.  
\elmma
\prf For $1\leq  i \leq k$, let $T_1^i,T_2^i,\cdots,T_{n-k+1}^i$ be  as above. 
 Take a nonzero vector $h_0 \in \clh_0^{k+1}$. By part (10) of proposition \ref{ppsn-spectrum}, one has
\begin{equation} \label{clh_0^{k+1}}
 \clh_0^k= \oplus_{\alpha_i \in \bbn}(T_1^k)^{\alpha_1}(T_2^k)^{\alpha_2} \cdots (T_{\ell_k-1}^{k})^{\alpha_{\ell_k-1}}\clh_0^{k+1}.
\end{equation}
Since  for $i\leq k$, $\clh_0^{k+1}\subset \clh_0^{i+1}$, there exists $t_i \in \bbbt$ 
such that 
 $T_{\ell_i}^ih_0=t_ih_0$ for all  $1 \leq i \leq k$.  By applying commutation relations (\ref{c3}), (\ref{c4}) and equation (\ref{<c_k}), 
 we observe that  $T_{\ell_i}^iT_j^k=T_j^kT_{\ell_i}^i$ on $\clh_0^k$.  Using this fact  and equation (\ref{clh_0^{k+1}}), one gets  $T_{\ell_i}^ih=t_ih$ for $h \in \clh_0^k$. Hence it follows from equation (\ref{decomposition}) that 
 if  $t^{'} \neq (t_1,t_2,\cdots, t_{k-1})$, $J(t^{'},a)$ is empty. 
 This  along with Lemma \ref{onea} proves that for only one value of 
 $(t,a) \in \bbbt^{k-1} \times \cla_{k-1}$, $J(t,a)$ is nonempty.
\qed \\
We call $t:=(t_1,t_2,\cdots t_k)$  and $\ell:=(\ell_1,\ell_2,\cdots ,\ell_k)$ defined as in above lemma the  angle and 
the rank respectively of the irreducible representation $\pi$ of $C^{n,k}$. 
The following lemma says that angle and rank completely determine an irreducible representation $\pi$ of $C^{n,k}$. 
\blmma \label{equivalence}  Assume that  the map $\Psi_m:C^{n,m} \rightarrow C(SU_q(n)/SU_q(n-m))$ sending $w_j^i$ to $u_j^i$  is an isomorphism for all $m=1,2,\cdots k-1$.
Let $\pi$ and $\pi^{'}$ be  irreducible representations of $C^{n,k}$ on a Hilbert space $\clh$ and $\clh^{'}$ respectively 
with same rank and same angle. Then $\pi$ and $\pi^{'}$ are equivalent.
\elmma
\prf Let $h$ and $h^{'}$ be nonzero vectors in $\clh_0^{k+1}$ and $(\clh_0^{k+1})^{'}$ respectively. Using remark \ref{ker} and part (9) and part (10) of proposition \ref{ppsn-spectrum}, 
we get   orthonormal bases for $\clh$ and $\clh^{'}$ given by 
\[
   \Big\{w_{\alpha}:=\frac{(T_1^1)^{\alpha_1^1}(T_2^1)^{\alpha_2^1}\cdots (T_{\ell_1-1}^1)^{\alpha_{\ell_1-1}^1}\cdots (T_1^k)^{\alpha_1^k}(T_2^k)^{\alpha_2^k}\cdots 
  (T_{\ell_k-1}^k)^{\alpha_{\ell_k-1}^k}h}{\|(T_1^1)^{\alpha_1^1}(T_2^1)^{\alpha_2^1}\cdots (T_{\ell_1-1}^1)^{\alpha_{\ell_1-1}^1}\cdots (T_1^k)^{\alpha_1^k}(T_2^k)^{\alpha_2^k}\cdots 
  (T_{\ell_k-1}^k)^{\alpha_{\ell_k-1}^k}h\|}\Big\}_{\alpha_i^j \in \bbn}
  \]
and 
\[
   \Big\{w_{\alpha}^{'}:=\frac{((T_1^1)^{'})^{\alpha_1^1}((T_2^1)^{'})^{\alpha_2^1}\cdots 
  ((T_{\ell_1-1}^1)^{'})^{\alpha_{\ell_1-1}^1}\cdots (T_1^k)^{'})^{\alpha_1^k}\cdots 
  ((T_{\ell_k-1}^k)^{'})^{\alpha_{\ell_k-1}^k}h^{'}}{\|((T_1^1)^{'})^{\alpha_1^1}(T_2^1)^{'})^{\alpha_2^1}\cdots 
  ((T_{\ell_1-1}^1)^{'})^{\alpha_{\ell_1-1}^1}\cdots((T_1^k)^{'})^{\alpha_1^k}\cdots 
  ((T_{\ell_k-1}^k)^{'})^{\alpha_{\ell_k-1}^k}h^{'}\|}\Big\}_{\alpha_k^l \in \bbn}
  \]
 respectively. Define 
 \begin{IEEEeqnarray}{rCl}
  U:\clh &\longrightarrow& \clh^{'} \nonumber \\
  w_{\alpha} &\longmapsto& w_{\alpha}^{'} \nonumber 
 \end{IEEEeqnarray}
 In the light of equation (\ref{decomposition}), 
 we only  need to show that for $1 \leq j \leq n$, 
 \[
  U\pi(w_j^{n-k+1})U^*=\pi^{'}(w_k^{n-k+1}).
 \]
 Using equation (\ref{c_k=0}) and  the fact that angle of $\pi$ and $\pi^{'}$ are same, one can show that 
\begin{equation} \label{e1}
 U\pi(w_j^{n-k+1})U^*h=\pi^{'}(w_j^{n-k+1})h^{'} \quad \mbox{  for all  } 1 \leq j \leq n.
\end{equation}
By applying commutation relations (\ref{c1}), (\ref{c2}), (\ref{c3}) and equation (\ref{e1}), we have $U\pi(w_n^{n-k+1})w_{\alpha}=\pi^{'}(w_n^{n-k+1})w_{\alpha}^{'}$ for all $\alpha$ and hence $U\pi(w_n^{n-k+1})U^*=\pi^{'}(w_n^{n-k+1})$. 
Assume the claim  for $j=n,n-1,\cdots i+1$. Again by applying  
commutation relation (\ref{c1}), (\ref{c2}), (\ref{c3}) and equation (\ref{e1}), we get $U\pi(w_i^{n-k+1})w_{\alpha}=\pi^{'}(w_i^{n-k+1})w_{\alpha}^{'}$ for all $\alpha$ and hence $U\pi(w_i^{n-k+1})U^*=\pi^{'}(w_i^{n-k+1})$. 
This completes the proof. 
\qed\\
Now we prove the main result of this paper.
\bthm
For $1\leq m \leq n$, the map  $\Psi_m:C^{n,m} \rightarrow C(SU_q(n)/SU_q(n-m))$ sending $w_j^i$ to $u_j^i$ for all $1\leq j \leq n$ and $n-m+1 \leq i \leq n$ is an isomorphism.
\ethm
\prf
We apply induction on $m$. For $m=1$, the claim is true thanks to proposition \ref{m=1}. We assume that the result is true for $m=1,2,\cdots,k-1$. 
We will show that for $m=k$, the map $\Psi_k$ is an  isomorphism. It is enough to show that $\Psi_k$ is injective. Take an irreducible representation $\pi$ 
of $C^{n,k}$.  Let $\ell=(\ell_1,\ell_2,\cdots,\ell_k)$ and $t=(t_1,t_1,\cdots,t_k)$ be the rank and the angle respectively  
of $\pi$. Define $a:=(n+1-\ell_1,n-\ell_2, \cdots ,n+2-k-\ell_k)$.  Consider the irreducible representation $\eta_{([t]_n,w(a))}\circ\Psi_k$ 
of the $C^*$-algebra $C^{n,k}$.
It is easy to check that rank and angle of this irreducible representation is $\ell$ and $t$ respectively. Hence by Lemma \ref{equivalence}, 
the representations $\pi$ and $\eta_{([t]_n,w(a))}\circ \Psi_k$  of $C^{n,k}$ are equivalent. 
This shows that any irreducible representation of $C^{n,k}$ factors through the map $\Psi_k$. Hence image of an element in $\ker(\Psi_k)$ 
under any irreducible representation of $C^{n,k}$ is zero which implies that the element is zero. Therefore, 
the map $\Psi_k$ is injective.
\qed \\
This  establishes  $C(SU_q(n)/SU_q(n-m))$ as 
the universal $C^*$-algebra given by  finite sets of generators and relations.

\noindent\begin{footnotesize}\textbf{Acknowledgement}:
I would like to thank  Arup Kumar Pal for making me work on the problem and for his valuable suggestions. Proof of 
proposition \ref{ppsn-spectrum} is taken from  his notes.
\end{footnotesize}

\noindent{\sc Bipul Saurabh} (\texttt{saurabhbipul2@gmail.com})\\
         {\footnotesize Institute of Mathematical Sciences, IV Cross Road, CIT Campus,
Taramani,
Chennai, 600113, INDIA}


\begin{thebibliography}{10}


\bibitem{ChaPal-2003aa}
Partha~Sarathi Chakraborty and Arupkumar Pal.
\newblock Equivariant spectral triples on the quantum {${\rm SU}(2)$} group.
\newblock {\em $K$-Theory}, 28(2):107--126, 2003.  

\bibitem{ChaPal-2003ac}
Partha~Sarathi Chakraborty and Arupkumar Pal.
\newblock Characterization of spectral triples: A combinatorial approach. arXiv:math.OA/0305157.
\newblock 2003.

\bibitem{ChaPal-2008aa}
Partha~Sarathi Chakraborty and Arupkumar Pal.
\newblock Characterization of {${\rm SU}_q(\ell+1)$}-equivariant spectral triples
  for the odd dimensional quantum spheres.
\newblock {\em J. Reine Angew. Math.}, 623:25--42, 2008.

\bibitem{KliSch-1997aa}
Anatoli Klimyk and Konrad Schm{\"u}dgen.
\newblock {\em Quantum groups and their representations}.
\newblock Texts and Monographs in Physics. Springer-Verlag, Berlin, 1997.



\bibitem{NesTus-2012ab}
Sergey Neshveyev and Lars Tuset.
\newblock Quantized algebras of functions on homogeneous spaces with {P}oisson
  stabilizers.
\newblock {\em Comm. Math. Phys.}, 312(1):223--250, 2012.



\bibitem{PodVai-1999aa}
G.~B. Podkolzin and L.~I. Vainerman.
\newblock Quantum {S}tiefel manifold and double cosets of quantum unitary
  group.
\newblock {\em Pacific J. Math.}, 188(1):179--199, 1999.

\bibitem{PalSun-2010aa}
Arupkumar Pal and S.~Sundar.
\newblock Regularity and dimension spectrum of the equivariant spectral triple
  for the odd-dimensional quantum spheres.
\newblock {\em J. Noncommut. Geom.}, 4(3):389--439, 2010.

\bibitem{Pod-1995aa}
Piotr Podle{{\'s}}.
\newblock Symmetries of quantum spaces. {S}ubgroups and quotient spaces of
  quantum {${\rm SU}(2)$} and {${\rm SO}(3)$} groups.
\newblock {\em Comm. Math. Phys.}, 170(1):1--20, 1995. 

\bibitem{Sau-2015aa}
Bipul Saurabh.
\newblock On quantum quaternion spheres, 2015. 
\newblock To appear in {\em Proc. Indian Acad. Sci. Math. Sci.}


\bibitem{Sau-2015ab}
Bipul Saurabh.
\newblock {$q$}-invariance of quantum quaternion spheres. 
\newblock arXiv:1510.01862, 2015.



\bibitem{VakSou-1988aa}
L.~L. Vaksman and Ya.~S. Soibelman.
\newblock An algebra of functions on the quantum group {${\rm SU}(2)$}.
\newblock {\em Funktsional. Anal. i Prilozhen.}, 22(3):1--14, 96, 1988.

\bibitem{VakSou-1990ab}
L.~L. Vaksman and Ya.~S. Soibelman.
\newblock Algebra of functions on the quantum group {${\rm SU}(n+1),$} and
  odd-dimensional quantum spheres.
\newblock {\em Algebra i Analiz}, 2(5):101--120, 1990.




\end{thebibliography}
\end{document}